\documentclass[10pt, reqno]{amsart}
\usepackage{amssymb,mathrsfs,amsmath}
 \usepackage{graphicx}
\usepackage{amscd,amsmath,amsopn,amssymb,amsthm,multicol}
\usepackage{mathtools}
\usepackage{xcolor}
\usepackage{mathrsfs}
\usepackage{color}
\usepackage[matrix, all, 2cell]{xy}
\usepackage{lscape}
\usepackage{slashed}
\usepackage{graphicx}
\usepackage{epsf}
\usepackage{setspace}
\usepackage{upgreek}
\usepackage{textgreek}
\usepackage{tikz}
\usepackage{multirow}
\usepackage{cancel}
\usepackage{soul}
\usepackage{comment}
\usepackage{wasysym}
\usepackage{bbm}
\usepackage{diagbox}
\usepackage{tikz-cd}
\usepackage{bbold}
\usepackage{booktabs}

\advance\oddsidemargin by -1.0cm
\advance\evensidemargin by -1.0cm
\advance\topmargin by -1.0cm

\numberwithin{equation}{section}

 \newcommand{\K}{\mathbb{K}}
 \newcommand{\R}{\mathbb{R}}
\newcommand{\Z}{\mathbb{Z}}
\newcommand{\C}{\mathbb{C}}
\newcommand{\Hn}{\mathbb{H}}

\renewcommand{\sp}{\mathfrak{sp}}

\DeclareMathOperator{\Hol}{\mathsf{Hol}}

\newcommand{\Gam}{\Gamma}

\DeclareMathOperator{\imm}{\mathsf{Im}}

\DeclareMathAlphabet{\mathscrbf}{OMS}{mdugm}{b}{n}

\DeclareMathOperator{\SO}{\mathsf{SO}}
\DeclareMathOperator{\Sp}{\mathsf{Sp}}
 \DeclareMathOperator{\SU}{\mathsf{SU}}

\DeclareMathOperator{\U}{\mathsf{U}}

\DeclareMathOperator{\G}{\mathsf{G}}
\DeclareMathOperator{\Lie}{\mathsf{Lie}}
\DeclareMathAlphabet{\mathpzc}{OT1}{pzc}{m}{it}
\DeclareMathOperator{\Hh}{\mathsf{H}}

\DeclareMathOperator{\E}{\mathsf{E}}
\DeclareMathOperator{\Gl}{\mathsf{GL}}
\DeclareMathOperator{\Sl}{\mathsf{SL}}

\DeclareMathOperator{\Aut}{\mathsf{Aut}}
\DeclareMathOperator{\id}{\mathsf{id}}
\DeclareMathOperator{\Ed}{\mathsf{End}}
\DeclareMathOperator{\Ad}{\mathsf{Ad}}
\DeclareMathOperator{\vol}{\mathsf{vol}}
\DeclareMathOperator{\Id}{\mathsf{Id}}
\newcommand{\fr}{\mathfrak}
\newcommand{\al}{\alpha}
\newcommand{\be}{\beta}
\newcommand{\mc}{\mathcal}

\newcommand{\ep}{\varepsilon}

\newcommand{\cc}{\big(}
\newcommand{\CC}{\Big(}
\newcommand{\rr}{\big)}
\newcommand{\RR}{\Big)}
\newcommand{\om}{\omega}
\newcommand{\Om}{\Omega}

\usepackage{amsmath}
\DeclareFontFamily{U}{mathx}{}
\DeclareFontShape{U}{mathx}{m}{n}{<-> mathx10}{}
\DeclareSymbolFont{mathx}{U}{mathx}{m}{n}
\DeclareMathAccent{\widehat}{0}{mathx}{"70}
\DeclareMathAccent{\widecheck}{0}{mathx}{"71}

\DeclareMathAlphabet{\mathscrbf}{OMS}{mdugm}{b}{n}

\newcommand{\Vg}{\ensuremath{\mathsf{V}}}
\newcommand{\Wg}{\ensuremath{\mathsf{W}}}
\newcommand{\Ug}{\ensuremath{\mathsf{U}}}

\newcommand{\Gg}{\ensuremath{\mathsf{G}}}
\newcommand{\Lg}{\ensuremath{\mathsf{L}}}
\newcommand{\Kg}{\ensuremath{\mathsf{K}}}
\newcommand{\Xg}{\ensuremath{\mathsf{X}}}
\newcommand{\Ss}{\ensuremath{\mathsf{S}}}

\DeclareMathOperator{\Ker}{\mathsf{Ker}}

\DeclareMathOperator{\dd}{d}

\newtheorem{theorem}{Theorem}[section]
\newtheorem{lem}[theorem]{Lemma}
\newtheorem{prop}[theorem]{Proposition}
\newtheorem{corol}[theorem]{Corollary}

\theoremstyle{definition}
\newtheorem{defi}[theorem]{Definition}
\newtheorem{example}[theorem]{Example}
 \newtheorem{rem}[theorem]{Remark}
 
\theoremstyle{remark}

\numberwithin{equation}{section}

\def\bd{\begin{defi}}
\def\ed{\end{defi}}
\def\bt{\begin{theorem}}
\def\et{\end{theorem}}
\def\bl{\begin{lem}}
\def\el{\end{lem}}
\def\bp{\begin{prop}}
\def\ep{\end{prop}}
\def\br{\begin{rem}}
\def\er{\end{rem}}
\def\bc{\begin{corol}}
\def\ec{\end{corol}}
\def\bex{\begin{example}}
\def\eex{\end{example}}
\def\pr{\begin{proof}}
\def\pro{\end{proof}}
\def\eqna{\begin{eqnarray*}}
\def\eqnaa{\begin{eqnarray}}
\def\deqna{\end{eqnarray*}}
\def\deqnaa{\end{eqnarray}}

 \definecolor{crew}{rgb}{0.2,0.5,0.2}
\definecolor{red}{rgb}{0.57,0.11,0.15}
\definecolor{cobalt}{rgb}{0.04,0.3,0.85}

\usepackage[backref=page]{hyperref}
\renewcommand*{\backref}[1]{}
\renewcommand*{\backrefalt}[4]{%
	\ifcase #1 %
		\or        (cited on page~#2)%
	\else      (cited on pages~#2) %
	\fi}
\hypersetup{colorlinks=true,linkcolor=cobalt, citecolor=blue,urlcolor=black}

\makeatletter
\def\subsubsection{\@startsection{subsubsection}{3}%
  \z@{.5\linespacing\@plus.7\linespacing}{.3\linespacing}%
  {\normalfont\bfseries}}
\makeatother

\begin{document}

\title[Submanifolds of almost quaternionic skew-Hermitian manifolds]{Submanifolds of almost quaternionic skew-Hermitian manifolds}
\author{Ioannis Chrysikos and Jan Gregorovi\v{c}}
\address{Department of Mathematics and Statistics, 
 	Faculty of Science, Masaryk University,	Kotl\'{a}\v{r}sk\'{a} 2, 611 37 Brno, 
Czech Republic. \\ ORCID ID: 0000-0002-0785-8021 }
\email{chrysikos@math.muni.cz}

\address{Department of Mathematics, Faculty of Science, University of Ostrava, 701 03 Ostrava, Czech Republic, and Institute of Discrete Mathematics and Geometry, TU Vienna, Wiedner Hauptstrasse 8-10/104, 1040 Vienna, Austria. \\ ORCID ID: 0000-0002-0715-7911}\email{jan.gregorovic@seznam.cz}

\keywords{Almost symplectic submanifolds; almost complex submanifolds; almost pseudo-Hermitian submanifolds; almost quaternionic submanifolds; quaternionic skew-Hermitian symmetric spaces}
\subjclass[2020]{53C05, 53C10, 53C15, 53C26, 53C30}

\begin{abstract}
We  investigate   several classes of submanifolds of almost quaternionic skew-Hermitian manifolds $(M^{4n}, Q, \om)$, including almost symplectic, 
almost complex,  almost pseudo-Hermitian and almost quaternionic submanifolds.
 In the torsion-free case, we  realize each type of submanifold considered  in the theoretical part  by constructing explicit examples of submanifolds of semisimple quaternionic skew-Hermitian symmetric spaces.
  \end{abstract}
 
 \maketitle 
 


\tableofcontents

\pagestyle{headings}






\section{Introduction and preliminaries}\label{sec1}

 \subsection{Motivation and summary}

Since its introduction by \'E. Cartan, the holonomy group of an affine connection has played a
key role in   differential geometry. 
On a symplectic manifold $(M^{2m}, \om)$  a torsion-free connection $\nabla$ preserving $\om$   is said to be
of \textsf{special symplectic holonomy}, if its holonomy group $\Hol(\nabla)$  is absolutely irreducible and  properly contained   in the symplectic group $\Sp(2m, \R)$.
Nowadays, it is well-known that special symplectic holonomy groups are examples of \textsf{exotic holomomies}, i.e.,  groups that were missing from   the classical list of Berger (see \cite{Br99}).  
Historically, such holonomy groups were the last to be discovered. The first examples  with special symplectic holomomy were constructed in dimension 4 (real or complex)  by  R. Bryant, by employing twistor methods \cite{Br91}.  The complete 
classification was  carried out later in a series of works by Q-S. Chi, S. Merkulov and L. Schwachh\"ofer (see  for example  \cite{CMS1, MS1} and a summary can be found in \cite[p.~3]{Schw}).  
By these works it is known that 
torsion-free connections with prescribed special symplectic holonomy admit 
 a universal construction, based on the realization of a certain  Poisson structure. 

The list of special symplectic holonomies includes, between other notable entries, the group $\Gg= \SO^*(2n)\Sp(1)$  (cf.~\cite[Table 5]{Schw2}),  where  $\SO^*(2n)$ denotes the quaternionic real form of  $\SO(2n, \C)$.\footnote{For clarity  we remark  that the group $\SO^*(2n)$ is  also denoted by $\SO(n, \Hn)$.} 
A comprehensive  study of such $\Gg$-structures was   given in \cite{CGWPartI, CGWPartII},  in  collaboration with H.~Winther,  where we  introduced  the notion of  \textsf{almost quaternionic skew-Hermitian manifolds} $(M^{4n}, Q, \om)$. These are $4n$-dimensional smooth manifolds with $n>1$ whose frame bundle admits a reduction to the Lie group $\SO^*(2n)\Sp(1)$. Such geometric structures are encoded by pairs $(Q, \om)$, where $Q\subset\Ed(TM)$ is an almost quaternionic structure on $M$ and $\om\in\Omega^{2}(M)$ is an almost symplectic form which is $Q$-Hermitian. Thus they form a ``symplectic analogue'' of the better understood almost quaternionic pseudo-Hermitian structures, i.e., $\Sp(p, q)\Sp(1)$-structures, see  also the diagram in Section \ref{section1}.\footnote{In this text  almost quaternionic skew-Hermitian manifolds are abbreviated as  \textsf{almost qs-H manifolds}.}    In the torsion-free case this yields  the notion of \textsf{quaternionic skew-Hermitian manifolds}, which are $4n$-dimensional manifolds endowed with a torsion-free connection with holonomy contained in $\SO^*(2n)\Sp(1)$ for $n>1$.

Among several results established in \cite{CGWPartI, CGWPartII}, we proved that  for $n>3$ there exist five  pure intrinsic torsion types $\mc{X}_1, \cdots, \mc{X}_5$ of $\SO^*(2n)\Sp(1)$-structures and  seven $\SO^*(2n)\Sp(1)$-irreducible intrinsic torsion components $\mc{X}_{1}, \ldots, \mc{X}_{7}$ of $\SO^*(2n)$-structures.  One  can then mix these irreducible components to control the (first-order)  integrability of the pairs $(Q, \om)$ and $(H, \om)$. 
For convenience, we summarize the most important types in a table, see Table \ref{Table1} below.

\vskip 0.5cm
\begin{table}[h!]
\caption{\bf Important types of almost qs-H manifolds/almost hs-H manifolds}
  \label{Table1}
\begin{tabular}{@{} l | l | l @{}}
\toprule
$(Q, \om)$ & $\text{intrinsic torsion type}$ & $\text{characterization}$ \\
\hline
&  $\mc{X}_{15}$ & $\dd\om=0$ (\text{symplectic}) \\
&  $\mc{X}_{345}$ & $Q$  \text{is quaternionic} \\
&  $\mc{X}_{5}$ & $\dd\om=0$ $\text{and}$  $Q$  \text{is quaternionic} \\
&  $\mc{X}_{145}$ & $\text{locally} \ \exists \ f \in C^{\infty}(M) \ \text{such that} \ f\om \ \text{is symplectic}$\\
&  $\mc{X}_{234}$ &  $\text{the intrinsic torsion is a 3-form}$\\
&  $\mc{X}_{4}$ &  $\text{the intrinsic torsion is of vectorial type}$\\
\hline
$(H, \om)$ &  $\mc{X}_{1567}$ & $\dd\om=0$ (\text{symplectic}) \\
&  $\mc{X}_{3457}$  &    $H$  \text{is hypercomplex} \\
&  $\mc{X}_{57}$ & $\dd\om=0$ $\text{and}$  $H$  \text{is hypercomplex} \\
&  $\mc{X}_{234}$ &  $\text{the intrinsic torsion is a 3-form}$\\
&  $\mc{X}_{47}$ &  $\text{the intrinsic torsion is of vectorial type}$\\
\bottomrule
\end{tabular}
\smallskip
\end{table}

\smallskip

From the viewpoint of holonomy theory  $\SO^*(2n)$-structures are  less important, as it is known  that    torsion-free affine connections with (irreducible) full holonomy group  
$\SO^*(2n)$ cannot exist, see  \cite{Br99}. However, such geometries  are still useful, as  they facilitate the understanding of $\SO^*(2n)\Sp(1)$-structures. 
More recently,  in  collaboration with V.~Cort\'es    \cite{CCG},  we studied
the second-order differential geometry  (curvature) of $\SO^*(2n)\Sp(1)$-structures, with an emphasis on the torsion-free case. We also    described a construction of $\SO(2(n+1))$-structures  of type $\mc{X}_{3457}$ on the Swann bundle associated to any   qs-H manifold $(M^{4n}, Q, \om)$.


In this work we  present  a systematic study  of various classes of submanifolds of almost qs-H manifolds $(M^{4n}, Q, \om)$, such as 
\textsf{almost symplectic submanifolds}, \textsf{almost complex submanifolds}, \textsf{almost pseudo-Hermitian submanifolds} and finally \textsf{almost quaternionic submanifolds}. 
In the torsion-free case part of our theory provides an analogue of certain results from the theory of submanifolds of quaternionic K\"ahler manifolds, see for example \cite{Gray, A68, F79,  BCU81,  Tasaki86, AM93, AM00a, AM002}.  

In particular, for a general almost qs-H manifold $(M^{4n}, Q, \om)$ we first develop a framework for studying almost symplectic submanifolds  of $(M^{4n}, Q, \om)$. 
These are even-dimensional  submanifolds $N^{2k}$ of $M^{4n}$ for which the induced 2-form $\hat\om:=\iota^*\om\in\Om^{2}(N)$ is assumed to be non-degenerate, where $\iota : N\to M$ denotes the  immersion. 
 Based on the fact that  $\om$ is an almost symplectic form,  in Section \ref{almost_symplectic_subQ1}  we consider the following   $\om$-orthogonal decomposition of $T M|_{N}$, 
\[
T M|_{N}=T N\oplus \nu(N)\,, 
\]
  where the  \textsf{normal bundle} $\nu(N)$ is defined by gluing together  the $\om$-orthogonal complement  $(T_{x}N)^{\perp_{\om}}$ of  $T_{x}N$ in $T_{x}M$, for all  $x\in N$.  This allows us to define the  \textsf{induced connection}  $\tilde\nabla_{X}Y:=(\nabla^{Q, \om}_{X}Y)^{\top}$,  where $\nabla^{Q, \om}$ is the unique minimal adapted connection on $(M, Q, \om)$ presented in \cite{CGWPartI}  (with respect to a ``normalization condition''),  and introduce the second fundamental form $\al$ by the Gauss formula $\nabla^{Q, \om}_{X}Y=\tilde\nabla_{X}Y+\al(X, Y)$, with $\al(X, Y)\in\Gamma(\nu(N))$ for any $X, Y\in\Gamma(TN)$. 
 We then show that $\tilde\nabla$  preserves $\hat\om$ and hence is  an almost symplectic connection on $N$. We further relate the torsion of $\tilde\nabla$ with the tangential part of the  torsion $T^{Q, \om}$ of $\nabla^{Q, \om}$. To obtain these results we first study  the extrinsic geometry of almost symplectic submanifolds of general almost Fedosov manifolds,  a topic of independent interest.  
 The details our outlined in Section \ref{Extrinsic_geom1}; see Propositions \ref{induced_connection} and \ref{second_fund_sym} which extend certain   results from \cite{CGRS09} about the extrinsic geometry of symplectic submanifolds of Fedosov manifolds, to the almost symplectic case.

Next we study  (compatible) \textsf{almost complex submanifolds} $(N^{2k}, \hat{J})$ of an almost qs-H manifold $(M^{4n}, Q, \om)$. By the term ``compatible'' we mean that the almost complex structure $\hat{J} : TN\to TN$ is locally    the restriction of an element  belonging to a (local) admissible basis $\{J_{a} : a=1, 2, 3\}$ of the almost quaternionic structure $Q$ (one may assume that $\hat{J}=J_{1}|_{TN}$ for example). Here we exploit the fact that $\nabla^{Q, \om}$ is an almost quaternionic connection, so locally there exist   1-forms $\gamma_{a}$ (depending on the choice of the admissible basis $\{J_a\}$), such that  $\nabla^{Q, \om}_{X}J_{a}=\gamma_{c}(X)J_{b}-\gamma_{b}(X)J_{c}$, 
  for any $X\in\Gamma(TM)$ and any cyclic permutation $(a, b, c)$ of $(1, 2, 3)$.
 Then we can obstruct the integrability of $\hat{J}$ in terms of the local 1-form 
 \[
 \uppsi:=(\gamma_3\circ J_1-\gamma_2)|_{TN}
 \]
  and the torsion $T^{Q, \om}$ of $\nabla^{Q, \om}$. This is the content of  Theorem \ref{Thm_1form}, which also treats the more general case of  almost qs-H manifolds of intrinsic torsion type $\mc{X}_{345}$ (i.e., when $Q$ is a quaternionic structure). In fact, when  $(M^{4n}, Q, \om)$ is either torsion-free  or of  type   $\mc{X}_{345}$, we  show that the obstruction to the integrability of the almost complex structure $\hat{J}$  lies in the $Q$-invariant  part $T^{Q}N$ of the tangent bundle  $TN$ of $N$ (see Proposition \ref{lemma_TQN}). In the torsion-free case this yields a similarity with the theory of almost complex submanifolds of quaternion K\"ahler manifolds with non-zero scalar curvature (cf. \cite{AM00a, AM002}).
 
 Given  an almost complex submanifold $(N^{2k}, \hat{J})$ of an almost qs-H manifold $(M^{4n}, Q, \om)$, the next natural step is to assume that the induced 2-form $\hat\om$ is non-degenerate  on $N$.  Under this assumption, Proposition \ref{almost_complex_sub}  shows that $N^{2k}$  admits   an almost pseudo-Hermitian structure   $(\hat{J}, \hat{g})$. 
   We then characterize the pair $(\hat{J}, \hat{g})$ in terms of the four pure Grey-Hervela classes of almost (pseudo) Hermitian structures (\cite{GrHerv}), using an equivariant projection of the torsion  $T^{Q, \om}$ and the local 1-form $\uppsi$. 
    This is the content of Corollary \ref{clas_corol_1},  Theorem \ref{Gr_H}, and Corollary \ref{Gr_H_corol}, where we explicitly characterize  the type of the  pair $(\hat{J}, \hat{g})$ (e.g., almost pseudo-K\"ahler, nearly pseudo-K\"ahler, pseudo-K\"ahler, etc).
          
    We also employ  the induced connection $\tilde\nabla$ to 
  obstruct the integrability of $\hat{J}$.  In  particular, for  an almost pseudo-Hermitian submanifold $(N^{2k}, \hat{J}, \hat{g})$ of an almost qs-H manifold $(M^{4n}, Q, \om)$,  we prove in Theorem   \ref{main_application_acs} that  $\hat{J}$ is integrable, provided that the following two conditions
    \[
   \pi_{J_1}(T^{Q, \om}(X, Y))=0\,,\quad \text{and} \quad
  \big(\gamma_3(X)J_2(Y)-\gamma_2(X)J_{3}(Y)\big)^{\top}=0
  \]
    hold for all $X, Y\in\Gamma(TN)$, where $\pi_{J_{1}}$ is  the projection  defined in Section \ref{admis_bas_int}.   When $(M^{4n}, Q, \om)$ is torsion-free we see that  $\hat{J}$ is integrable if and only if the second condition is satisfied. In Theorem \ref{totally_complex} we finally   establish a precise relationship between  this condition and the following ones:
\[
 \gamma_{2}|_{T_xN}=\gamma_{3}|_{T_xN}=0\,,\qquad  
 J_{2}T_{x}N\perp_{\om} T_{x}N\,.
\]
Here, the second condition resembles  the usual  \textsf{totally-complex} condition from the theory of    almost Hermitian submanifolds.

 In Section \ref{almostqsub} we  analyze \textsf{almost quaternionic submanifolds}   of torsion-free almost qs-H manifolds, i.e.,  qs-H manifolds $(M^{4n}, Q, \om)$. 
These are  submanifolds $N\subset M^{4n}$ whose tangent space $T_{x}N$ is $Q$-invariant at any $x\in N$, and hence they are  necessarily $4k$-dimensional for some $k<n$.  As a consequence of the general theory of almost quaternionic submanifolds of quaternionic manifolds (see \cite{AM93,Pet}), it follows that any such submanifold $N^{4k}\subset M^{4n}$
 is  a quaternionic manifold itself,  with the quaternionic structure defined by restriction, i.e., $\hat{Q}:=Q|_{TN}$. Moreover, such a submanifold  is totally geodesic with respect to any quaternionic connection on $M$. 
In our setting, we prove that any almost quaternionic submanifold  $(N^{4k}, \hat{Q}=Q|_{TN})$  of a qs-H manifold $(M^{4n}, Q, \om)$  for which the induced 2-form $\hat\om$ is non-degenerate  is itself  a quaternionic skew-Hermitian manifold. In other words,  $N^{4k}$ admits a torsion-free $\SO^*(2k)\Sp(1)$-structure $(\hat{Q}, \hat{\om})$, as specified in the Theorem \ref{theorem1}. 
    
In the final section, we   illustrate  the four main classes of submanifolds introduced above in the torsion-free case.  This is done by constructing homogeneous examples based on simple quaternionic skew-Hermitian symmetric spaces (see \cite{CGWPartI, CCG}). These examples are presented in  Proposition  \ref{ExThm1},   Examples \ref{Lan1} and \ref{Ex_SOstar} (note that   Example \ref{Lan1}  highlights a Lagrangian example).  
    
    The article is structured as follows: To set the stage we start with necessary preliminaries about $\SO^*(2n)\Sp(1)$-structures and $\SO^*(2n)$-structures (see Sections \ref{section1} and \ref{adapted_conn1}). In Section \ref{admis_bas_int}  we present a  result regarding the integrability of a local almost complex structure $J_{a}$ in an admissible basis of the almost quaternionic structure $Q$ (see Theorem \ref{them_1.7} and Corollary \ref{corol_Ja_integrable}). This will play an important role  in Section \ref{section2},  which is devoted to submanifolds. 
    Section \ref{section4} contains  our examples.

\smallskip
\noindent{\bf{Acknowledgements}}
The first author acknowledges  support  from the Czech Grant Agency (project GA24-10887S), and the  Horizon Europe MSCA project CaLIGOLA,   ID 101086123. 
He also thanks the Institute of Mathematics (UZH) at the University of Z\"urich for its hospitality during his academic stay in the Winter  semester 2025/2026. 
The second author acknowledges the support by the Austrian Science Fund (FWF) 10.55776/PAT4819724.

\subsection{A review of $\SO^*(2n)$-structures and  $\SO^*(2n)\Sp(1)$-structures}\label{prelim_1}\label{section1}
The  Lie group $\SO^*(2n)$  is defined as the  intersection 
\[
\Gl(n, \Hn)\cap\SO(2n, \C)\,,\quad\text{or}\quad \Gl(n, \Hn)\cap\Sp(4n, \R)
\]
 and
for $n=1$  it coincides with the compact Lie group $\SO(2)=\U(1)$.   From now on we will assume that  $n>1$ and view $\SO^*(2n)$  as a  non-compact real form of $\SO(2n, \C)$.
The  Lie  group $\SO^*(2n)$ has (real) dimension $n(2n-1)$ and is referred to as the \textsf{quaternionic real form} of  $\SO(2n, \C)$. We also set  
\[
\SO^*(2n)\Sp(1):=\SO^{\ast}(2n)\times_{\Z_2}\Sp(1)=(\SO^{\ast}(2n)\times\Sp(1))/{\Z_2}\,.
\]
 
 Roughly speaking, the $\Gg$-structures associated with the above two groups  provide  a symplectic analogue of the better understood almost hypercomplex  pseudo-Hermitian structures and almost quaternionic pseudo-Hermitian structures, which are $\G$-structures for the Lie groups $\Sp(p, q)$ and $\Sp(p,q)\Sp(1)$, respectively (see for example  \cite{Salamon86, AM, AC} for details on these  $\Gg$-structures). 
 
In order to facilitate understanding, we summarize all these  geometric structures  in a diagram, together with the   inclusion relations between them, see below.   In this diagram  the subscript 
$0$ indicates the defining tensors of the   linear model associated to  the  depicted $\Gg$-structure.
On the other hand, each inclusion refers to the appropriate choice of the
structures.  For example, the inclusion $\Gl(n, \Hn)\hookrightarrow \Gl(n, \Hn)\Sp(1)$ refers to  the almost quaternionic structure  
generated by an almost hypercomplex structure,  and more details can be found in   \cite{AM}. 

 In what follows, we will  restrict our attention to the right-hand-side of this diagram, which encodes  the skew-Hermitian   classes, with an emphasis on the almost qs-H case.
For numerous details omitted below, we refer the reader to  \cite{CGWPartI, CCG}, see also \cite[Section 2]{AM}  for standard facts  on almost hypercomplex and almost quaternionic manifolds  that we  avoid to recall.


{\small
\[
\scalebox{0.85}{
\xymatrix@C=0.6cm@R=0.8cm{
 &  \overset{\Gl(n,\Hn)\Sp(1)=\Aut(Q_0)}{\text{Almost Quaternionic}} & \\
 &  \overset{\Sl(n,\Hn)\Sp(1)=\Aut(Q_0, \vol_0)}{\text{Almost Unimodular Quaternionic}} \ar[u] & \\
 \overset{\Sp(p,q)\Sp(1)=\Aut(Q_0, g_0)}{\text{Almost Quaternionic Hermitian}} \ar@/^4pt/[ur] 
 & & \overset{\SO^{*}(2n)\Sp(1)=\Aut(Q_0, \omega_0)}{\text{Almost Quaternionic Skew-Hermitian}} \ar@/_4pt/[ul]\\
 &  \overset{\Gl(n,\Hn)=\Aut(H_0)}{\text{Almost Hypercomplex}} 
     \ar@/^85pt/@{->}[uuu]|{\color{red}} & \\
 &  \overset{\Sl(n,\Hn)=\Aut(H_0, \vol_0)}{\text{Almost Unimodular Hypercomplex}}  
     \ar[u]  
     \ar@/_60pt/@{->}[uuu]|{\color{blue}} & \\
 \overset{\Sp(p,q)=\Aut(H_0, g_0)}{\text{Almost Hypercomplex Hermitian}} \ar@/^4pt/[ur] 
  \ar@/^10pt/@{->}[uuu] & & \overset{\SO^{*}(2n)=\Aut(H_0, \omega_0)}{\text{Almost Hypercomplex Skew-Hermitian}} \ar@/_4pt/[ul] 
  \ar@/_10pt/@{->}[uuu]
}}
\]
}

  Let  us fix   a real vector space  $\Vg$  of dimension $4n$. As we mentioned above,  for our goals it suffices to assume that  $n>1$. Recall that  a \textsf{linear hypercomplex structure} on $\Vg$ is a triple  $H=\{J_1, J_2, J_3\}$ of anti-commuting linear complex structures $J_i\in\Ed(\Vg)$   with  $   J_1\circ J_2=J_3=-J_2\circ J_1$. 
On the other hand,   a \textsf{linear  quaternionic structure} on $\Vg$  is a 3-dimensional subalgebra   $Q \subset\Ed(\Vg)$ of the Lie algebra $\Ed(\Vg)$  generated by a
 linear hypercomplex structure $H$ on $\Vg$, i.e.,  $Q=\langle H\rangle$.  Obviously,  $Q$ is a Lie subalgebra  isomorphic  to  $\fr{sp}(1)\cong{\rm Im}(\Hn)$, where  as usual $\Hn$ denotes the  (division) algebra of quaternions.  
\bd  
On a hypercomplex vector space $(\Vg, H)$ (respectively, quaternionic vector space $(\Vg, Q=\langle H\rangle)$) a  non-degenerate  2-form $\omega \in\bigwedge^2\Vg^*$  which is $H$-Hermitian  (respectively, $Q$-Hermitian) is referred to as a \textsf{scalar 2-form} with respect to $H$ (respectively, with respect to $Q$).
\ed
 Consider the Lie group $\SO(2n, \C)$ of complex linear transformations preserving  
 the standard complex Euclidean metric on $\E:=\C^{2n}$ and identify $\E$ with  the 
 standard representation of $\SO^{\ast}(2n)\subset \SO(2n, \C)$. 
 Set also  $\Hh:=\C^2$ for the standard representation of $\Sp(1)$. Both  $\E, \Hh$   
 are of quaternionic type and the  standard  representation of $\SO^{\ast}(2n)\Sp(1)$ is 
   the real form   $[\E\Hh]$  inside $\E\otimes_{\C}{\Hh}$, fixed by the real structure $\epsilon_{\E}\otimes \epsilon_{\Hh}$, where    $\epsilon_{\E}, \epsilon_{\Hh}$    are the corresponding complex anti-linear involutions.  
There is a  standard linear quaternionic structure on $[\E\Hh]$, denoted by  $Q_0\cong\fr{sp}(1)$.
 Let $\om_0$ be the non-degenerate $Q_0$-Hermitian  linear 2-form on $[\E\Hh]$, given by
 $\om_0=g_{\E}\otimes\om_{\Hh}$, where $\om_{\Hh}$ is the generator of $[\bigwedge^2\Hh]^*$ and $g_{\E}$ is the generator of a 1-dimensional $\SO^*(2n)$-submodule in $[S^2\E]^*$, see  \cite[Prop.~2.6]{CGWPartI}.
 Then $([\E\Hh], Q_0, \omega_0)$ is a quaternionic skew-Hermitian vector space.   
 The 2-form $\om_0$ is called the  \textsf{standard scalar 2-form} on $[\E\Hh]$,  and we can define pseudo-Euclidean metrics
  $g_0^{a}$ by $g_0^{a}:= \om_0(\cdot , \mc{J}_a)$, 
where  $H_{0}:=\{\mc{J}_{a} : a=1, 2, 3 \}$ is the \textsf{standard admissible basis} of $Q_0$.   The triple $([\E\Hh], Q_0, \om_0)$ serves as the \textsf{linear model} of $\SO^*(2n)\Sp(1)$-structures, for more details see   \cite[Section 2]{CGWPartI}.

  \bd
{\rm(a)} Let $(M^{4n}, H)$ be an almost hypercomplex manifold.  A \textsf{scalar 2-form with respect to $H$} is a non-degenerate smooth 2-form $\om\in\Om^2(M)=\Gamma(\bigwedge^2T^*M)$  which is $H$-Hermitian, i.e.,
 \begin{equation}\label{QHerm}
\omega(J_{a}X, J_{a}Y)=\omega(X, Y)\,,\quad a=1, 2, 3\,.
\end{equation}
{\rm(b)} Let $(M^{4n}, Q)$ be an almost quaternionic manifold.  A \textsf{scalar 2-form with respect to $Q$} is a non-degenerate smooth 2-form $\om\in\Om^2(M)$  which is $Q$-Hermitian. This means that  $\om$  satisfies (\ref{QHerm})  for any (local) admissible basis $\{J_a : a=1, 2, 3\}$ of $Q$ and   $X, Y\in\Gamma(TM)$. \ed
 
\bd Let $M$ be a $4n$-dimensional manifold with $n>1$.\\
(a) An \textsf{almost hypercomplex skew-Hermitian  structure}  (in short, \textsf{almost hs-H structure}) on $M$ is reduction of the frame bundle $\mc{F}(M)$ of $M$ to $\SO^*(2n)$, that is, a $\SO^*(2n)$-structure. This is equivalent to say that 
$M$ admits a pair $(H, \omega)$, where $H=\{J_a : a=1, 2, 3\}$  is an almost hypercomplex structure on $M$ and $\omega$ is a scalar 2-form with respect to $H$.	In this case   $(M, H, \om)$  is said to be an \textsf{almost hypercomplex skew-Hermitian manifold} (in short, \textsf{almost hs-H manifold}). \\
(b)  An \textsf{almost quaternionic skew-Hermitian    structure} (in short, \textsf{almost qs-H structure}) on $M$ is a $\SO^*(2n)\Sp(1)$-structure.  This is equivalent to say that  $M$ admits a pair $(Q , \omega)$, where $Q\subset\Ed(TM)$ is an almost quaternionic structure  on $M$ and $\omega$ is  scalar 2-form with respect to $Q$.  In this case   $(M, Q, \om)$ is said to be an \textsf{almost quaternionic skew-Hermitian manifold} (in short, \textsf{almost qs-H manifold}).\\
(c) A \textsf{quaternionic skew-Hermitian transformation} between two almost qs-H manifolds $(M_i, Q_i, \om_i)$ $(i=1, 2)$ is a diffeomorphism
$f : M_1\to M_2$ such that $f^*\omega_2=\om_1$ and $f^*Q_2=Q_1$.  A \textsf{hypercomplex skew-Hermitian transformation} between two almost hs-H manifolds $(M_i, H_i, \om_i)$ $(i=1, 2)$ is defined in a similar manner.   If the scalar 2-forms  $\omega_i$ are both symplectic, then a qs-H  transformation  $f : (M_1, Q_1, \om_1)\to (M_2, Q_2, \om_2)$ is said to be  a  \textsf{quaternionic symplectomorphism}.
 \ed
\br
The terminology for $\SO^*(2n)$- and $\SO^*(2n)\Sp(1)$-structures  follows our previous works and is motivated  by Harvey's   classical description of the eight types of ``inner product spaces'' \cite{H90} (see also Remark \ref{terminology_1} below).
\er

 Fix an almost quaternionic skew-Hermitian manifold $(M, Q, \om)$.  According to  \cite[Proposition 3.9]{CGWPartI}, at any point $x\in M$  we may  identify $T_{x}M\cong[\E\Hh]$. 
The pair $(Q, \om)$ induces a volume form on $M$, namely $\vol:=\om^{2n}$.  Hence any almost qs-H manifold is oriented and in particular, $(Q, \om)$  gives rise to a so-called  \textsf{almost quaternionic unimodular structure}, which is encoded by the pair  $(Q, \vol=\om^{2n})$. This explains the arrows in the right-hand-side of the diagram given above. 
Moreover, choosing an admissible  basis $H=\{J_{a} : a=1, 2, 3\}$ of $Q$ we  can introduce three pseudo-Riemannian metrics of signature $(2n, 2n)$ on $M$,  defined by 
\[
g_{J_a}(X, Y):=(\om\circ J_{a})(X, Y)=\omega(X,  J_{a}Y)\,,\quad a=1, 2, 3
\]
 Next we will usually  write $g_{a}$, instead of $g_{J_{a}}$. 
Typically,  the metrics $g_a$  are only {\it locally} defined, since the admissible basis $H=\{J_a\}$ is only locally defined.  They are global tensors when $H$ is an almost hypercomplex structure, that is, when  $(H, \omega)$  is an    almost hypercomplex skew-Hermitian  structure on $M$.   Note that  the metrics $g_{a}$ are {\it not} $H$-Hermitian, although each $g_{a}$ is $J_a$-Hermitian, that is,  $g_{a}(J_{a}X, J_{a}Y)=g_a(X, Y)$ for all $a=1, 2, 3$ and $X, Y\in\Gamma(TM)$. 

Although the individual metrics $g_a$ are   only locally defined,  we mention that  any almost quaternionic skew-Hermitian manifold $(M, Q, \om)$ is endowed with some important   globally defined tensor fields.  However, we shall not use these tensors in the present work (interested readers may consult  \cite{CGWPartI, CCG} for further details).

 \subsection{Adapted connections and integrability conditions}\label{adapted_conn1}
 Let $(M, H=\{J_{a} : a=1, 2, 3\})$  be an  almost hypercomplex manifold.  An  \textsf{almost hypercomplex connection} on $M$   is a connection $\nabla$ on the tangent bundle $TM$ of $M$ which preserves $H$, that is, $\nabla J_{a}=0$ for all  $a=1, 2, 3$.  An almost hypercomplex connection  $\nabla$ which  is torsion-free  is called a  \textsf{hypercomplex connection}.
It is known that  there is a unique  almost hypercomplex connection $\nabla^H$ adapted to an almost hypercomplex structure  $H=\{J_a : a=1, 2, 3\}$, referred to as the \textsf{Obata connection},  whose torsion tensor is given by  $T^{H}=\frac{1}{12}\sum_{a=1}^{3}[[J_a, J_a]]$. 
 Here, $[[J_a, J_a]]$
denotes the Nijenhuis bracket of $J_a\in H$,   with $[[J_{a}, J_{a}]]=2N_{J_{a}}$, where $N_{J_{a}}$ is the Nijenhuis tensors of $J_a\in H$. 
Obata proved that   $H$ is \textsf{hypercomplex} if and 
 if and only if $\nabla^H$ is a hypercomplex connection, i.e.,  $T^H=0$. If  this is the case, then $(M, H)$  is said to be a \textsf{hypercomplex manifold}, see also \cite{AM}.

  Let $(M, Q)$  be an almost quaternionic manifold.
Clearly, any connection $\nabla$ on $TM$ induces a covariant derivative on all tensor bundles build from $TM$ (and its dual), and in particular on $\Ed(TM)\cong TM\otimes T^*M$.  We maintain the same notation for this connection, defined by
\[
(\nabla_{X}\sigma)Y=\nabla_{X}(\sigma Y)-\sigma(\nabla_{X}Y)\,, 
\]
for any $X, Y\in\fr{X}(M)$ and $\sigma\in\Gamma(\Ed(TM))$. 
Saying that $\nabla$ preserves the subbundle $Q\subset\Ed(TM)$ we mean that maps sections of $Q$ to sections of $T^*M\otimes Q\subset T^*M\otimes\Ed(TM)$.  Hence for any section $\sigma\in\Gamma(Q)$ we require that $\nabla\sigma$ is a 1-form with values in $Q\subset\Ed(TM)$, i.e., $\nabla\sigma\in\Gamma(T^*M\otimes Q)$ for all $\sigma\in\Gamma(Q)$. 
 Such a connection $\nabla$ is called an  \textsf{almost quaternionic connection} (subordinated to the quaternionic structure $Q$). For $n>1$ an almost quaternionic connection which  is torsion-free   is said to be a  \textsf{quaternionic connection} and in this case $(M. Q)$ is called a    \textsf{quaternionic manifold} (see  \cite{AM} for the special case $n=1$).  Note that the non-uniqueness of quaternionic connections is a standard result in quaternionic geometry,  see for example \cite{O84, Salamon86, AM}. 
 
 For the case of   \textsf{unimodular almost quaternionic structures} $(Q, \vol)$,   there exists a unique almost quaternionic connection $\nabla^{Q, \vol}$ preserving  the volume form $\vol$, whose torsion, locally in the domain of any admissible basis $\{J_a\}$ of $Q$, is given 
by 
\[
T^Q=T^H+\sum_{a=1}^{3}\partial(\tau_{a}^{H}\otimes J_{a})\,.
\]
Here, $\tau_{a}^{H}(X):=\frac{1}{4n-2}{\rm tr}(J_{a}T^{H}_{X})$ and $\partial$  denotes the operator of alternation. 
 This is  referred  to as the \textsf{unimodular Oproiu connection} and it can be proved  that $Q$ is quaternionic if and only  $T^{Q}=0$ (see \cite{AM}).  
 Note that $T^Q$ does not depend on the volume form. 

  Let us finally recall useful facts concerning \textsf{almost symplectic connections}.
Let $(M, \om)$ be an almost symplectic manifold. As before, we will maintain assuming that $\dim M>4$.  An affine connection $\nabla$ on $M$ which preserves $\om$, i.e., $\nabla\om=0$, is called an \textsf{almost symplectic connection}.
When in addition  $\nabla$ is torsion-free,  $T^{\nabla}=0$ identically, then $\nabla$ is referred to as a \textsf{symplectic connection}.
It is easy to see that for the torsion-free case, the condition $\nabla\om=0$ implies that $\dd\om=0$;     hence in this case $(M, \om)$ is a symplectic manifold.    
The converse is also true, see \cite[Proposition 2.5]{HHR07}.

Let us now introduce the following definition, which extends the classical definition of the so-called \textsf{Fedosov manifolds} (\cite{GRS98}).

\bd
An \textsf{almost Fedosov manifold} is an almost symplectic manifold $(M, \om)$ equipped with  an almost symplectic connection $\nabla$. 
When $\nabla$ is torsion-free, then  $(M, \om)$  is said to be  a \textsf{Fedosov manifold},  usually  denoted by $(M, \om, \nabla)$.
\ed
In other words, a Fedosov manifold is a symplectic manifold $(M, \om)$ endowed with a symplectic connection $\nabla$.  
Note that almost symplectic connections always exist, but in general are not unique. In particular,  almost symplectic connections with fixed  torsion are parametrized by symmetric 3-tensors, see for example \cite{Ton61, V86, HHR07}.

  Let us now return to $\SO^*(2n)$- and $\SO^*(2n)\Sp(1)$-structures  and recall   the associated    adapted connections,  denoted by $\nabla^{H, \omega}$ and $\nabla^{Q, \omega}$, respectively.  We mention that
 \begin{itemize} 
  \item For an almost hypercomplex skew-Hermitian manifold $(M, H, \om)$, the connection $\nabla^{H, \om}$ is both almost symplectic and   almost hypercomplex. It is referred to as an \textsf{almost hypercomplex skew-Hermitian connection} (in short,  \textsf{almost hs-H connection}).
  \item For an almost quaternionic skew-Hermitian manifold $(M, Q, \om)$, the connection $\nabla^{Q, \om}$ is both almost symplectic and almost quaternionic.
  It is referred to as an \textsf{almost quatenionic skew-Hermitian connection} (in short, \textsf{almost hs-H connection}).
  \end{itemize}
 In  \cite[Theorem 5.3]{CGWPartI} it was shown that both  $\nabla^{H, \om}, \nabla^{Q, \om}$ are \textsf{minimal adapted connections} for the corresponding $\Gg$-structures,  unique with respect to certain \textsf{normalization conditions}.  Let us recall, following \cite{CGWPartI},  their expression in the general case and  highlight   the torsion-free case according to the results in \cite{CGWPartII}.
\bt \label{connections} \textnormal{(\cite{CGWPartI, CGWPartII})}
{\rm (1)}  On an  almost hypercomplex skew-Hermitian manifold  $(M, H, \om)$  the  adapted connection '$\nabla^{H, \om}$ is given by $\nabla^{H, \om}=\nabla^{H}+A$,
where $\nabla^{H}$ is the  Obata connection associated to $H$, and $A$ is the $(1, 2)$-tensor field on $M$ defined by  
\[
\omega\cc A(X, Y), Z\rr=\frac{1}{2}(\nabla^{H}_{X}\omega)(Y, Z)\,,\quad X, Y, Z\in\Gamma(TM)\,.
\]
The connection $\nabla^{H, \omega}$ is torsion-free if and only if  
 \begin{equation}\label{omegTH}
 T^{H}=0\,,\quad \text{and}\quad  \nabla^{H}\omega=0\,.
 \end{equation}
In other words, $T^{H, \omega}=0$ if and only if $H$ is  hypercomplex and $\nabla^{H}$ is a symplectic connection. Moreover, when $T^{H, \om}=0$, then the scalar 2-form  $\om$ is a symplectic form.\\
\noindent{\rm (2)} On an almost quaternionic skew-Hermitian manifold  $(M, Q, \om)$  the adapted connection $\nabla^{Q, \om}$ is given by $\nabla^{Q, \om}=\nabla^{Q, \vol}+A$,
where $\nabla^{Q, \vol}$ is the    unimodular  Oproiu  connection associated to the pair $(Q, \vol=\om^{n})$, and $A$ is the $(1, 2)$-tensor field on $M$ defined by the relation
\[
\omega\cc A(X, Y), Z\rr=\frac{1}{2}(\nabla^{Q, \vol}_{X}\omega)(Y, Z)\,,\quad X, Y, Z\in\Gamma(TM)\,.
\]
The connection $\nabla^{Q, \omega}$ is torsion-free if and only if
 \begin{equation}\label{omegTHH}
 T^{Q}=0\,,\quad \text{and}\quad  \nabla^{Q, \vol}\omega=0\,.
 \end{equation}
 In other words, $T^{Q, \omega}=0$ if and only if $Q$ is quaternionic and $\nabla^{Q, \vol}$ is a symplectic connection.  
Moreover,  when  $T^{Q, \om}=0$, then  the scalar 2-form $\om$ is a symplectic  form.
\et

\br \label{terminology_1}
   As is customary in the torsion-free case,   the prefix ``almost'' can be omitted.   Thus, in this case one may simply refer to \textsf{hypercomplex skew-Hermitian structures} and \textsf{quaternionic skew-Hermitian structures}, a terminology also adopted in \cite{CCG}.
\er

\subsection{On the integrability of an element in an admissible basis}\label{admis_bas_int}
In preparation for the analysis in the forthcoming  section, 
 it will be useful to derive obstructions for the integrability of a member $J_{a}$ in an admissible basis of the almost quaternionic structure $Q$.
Note that such obstructions can be of  independent interest.

 Let $(Q, \om)$ be an almost qs-H structure and let $H=\{J_a : a=1, 2, 3\}$ be an  admissible basis of $Q$.  Since the almost qs-H connection $\nabla^{Q, \om}$ preserves $Q$, there are local 1-forms $\gamma_{a}$, depending on the choice of the admissible basis $\{J_a\}$, such that 
\begin{equation}\label{quat_conn}
\nabla^{Q, \om}_{X}J_{a}=\gamma_{c}(X)J_{b}-\gamma_{b}(X)J_{c}
\end{equation}
  for any $X\in\Gamma(TM)$ and any cyclic permutation $(a,b,c)$ of $(1, 2, 3)$.
     We also consider  the projections   $\pi_{J_{a}}$,  
defined by
\[
\pi_{J_{a}}(\phi)(X,  Y):=\frac{1}{4}\CC\phi(X,  Y)+J_{a}\big(\phi(J_{a}X, Y)+\phi(X,  J_{a}Y)\big)-\phi(J_{a}X, J_{a}Y)\RR\,,\quad a=1, 2, 3,
\]
 for any $X,  Y\in  [\E\Hh]$ and $\phi\in\bigwedge^2[\E\Hh]^*\otimes [\E\Hh]$. Furthermore, we set  
\[
\pi_{H}:=\frac{2}{3}(\pi_{J_{1}}+\pi_{J_{2}}+\pi_{J_{3}})=\frac{2}{3}\sum_{a=1}^{3}\pi_{J_{a}}\,.
\]
This operator plays a fundamental role in almost hypercomplex geometry, where it is used to decompose the corresponding intrinsic torsion into its irreducible components, see \cite{G97}.
It is also important in the setting of $\SO^*(2n)\Sp(1)$-structures, see \cite[pp.~2648--2649]{CGWPartI} for details. 

 \bl\label{integrable_J_a}
Let $(M, Q, \om)$ be an almost qs-H manifold and let $H=\{J_a : a=1, 2, 3\}$ be an admissible basis of $Q$.
Then the local almost complex structure $J_a\in H$ is integrable if and only if $\pi_{J_{a}}(T^{H, \om})=0$.
 \el
 \pr
  Of course, the choice of the admissible basis $H=\{J_a : a=1, 2, 3\}$ provides locally the adapted connection $\nabla^{H,\omega}$ of the  locally defined almost hs-H structure $(H=\{J_{a} : a=1, 2, 3\},\om)$. 
By definition,  $\nabla^{H, \om}$ preserves $H$, and hence we have $\nabla^{H, \om}J_{a}=0$. 
Therefore,  for the Nijenhuis tensor $N_{J_{a}}$  of $J_{a}$ we see that 
 \begin{eqnarray*}
N_{J_{a}}(X,Y)&=&(\nabla^{H, \om}_{J_{a}X} J_{a})Y - (\nabla^{H, \om}_{J_{a}Y} J_{a})X
  + J_{a}\big((\nabla^{H, \om}_Y J_{a})X - (\nabla^{H, \om}_X J_{a})Y\big)  \\
&& - T^{H, \om}(J_{a}X, J_{a}Y) + J_{a}\,T^{H, \om}(J_{a}X, Y) + J_{a}\,T^{H, \om}(X, J_{a}Y) + T^{H, \om}(X, Y)\,,\\
&=&4\pi_{J_{a}}(T^{H, \om})(X, Y)\,.
\end{eqnarray*}
 Thus $N_{J_{a}}(X, Y)=4\pi_{J_{a}}(T^{H, \om})(X, Y)$ for all  $X, Y\in\Gamma(TM)$,  and  our claim follows.
 \pro

  Let us now consider  the difference $\nabla^{H, \om}-\nabla^{Q, \om}$ and also 
  the difference $T^{H, \om}-T^{Q, \om}$  between the torsion forms $T^{H, \om}$ and $T^{Q,  \om}$ of $\nabla^{H, \om}$ and $\nabla^{Q,\om}$, respectively. 
These differences are isomorphic to the torsion in the intrinsic torsion component $\mc{X}_{67}$ of $\SO^*(2n)$-structures. 
 Moreover, by part (1) in  \cite[Theorem 5.3]{CGWPartI} we   know that $\mc{X}_6\subset\imm(\pi_{H})$, while $\mc{X}_7\subset {\sf Ker}(\pi_{H})$. Note that, in general,  we have that $\imm(\pi_{H})=\mc{X}_{126}$ and $\mc{X}_{3457}\subset{\sf Ker}(\pi_{H})$.

  \bt\label{them_1.7}
 Let $H=\{J_a\}$ be an admissible basis of $Q$ and  let $\gamma_{a}$ $(a=1, 2, 3)$ be  local 1-forms such that the relation {\rm(\ref{quat_conn})} is satisfied. Then,  the adapted connections $\nabla^{H,\omega}$ and $\nabla^{Q, \om}$ are related by
\begin{equation}\label{1.7}
\nabla^{H, \om}_X=\nabla^{Q, \om}_X-\frac12\sum^3_{a=1}\gamma_a(X)J_a-\frac23\pi_S\big(\omega\otimes (\sum^3_{a=1}\gamma_a(J_a))^T\big)(X)\,,
\end{equation}
where  $T$ denotes the symplectic transpose and $\pi_{S}$ is the projection introduced in \cite[p.~2641]{CGWPartI}.
The torsion forms $T^{H, \om}$ and $T^{Q, \om}$ satisfy 
 \begin{align*}
T^{H, \om}(X,Y)&=T^{Q, \om}(X,Y)+\frac12\sum^3_{a=1}(\gamma_a(Y)J_aX-\gamma_a(X)J_aY)\\
&-\frac23\left\{\pi_S\big(\omega\otimes (\sum^3_{a=1}\gamma_a(J_a))^T\big)(X,Y)-\pi_S\big(\omega\otimes (\sum^3_{a=1}\gamma_a(J_a))^T\big)(Y,X)\right\}
\end{align*}
where by $\gamma_{a}(J_a):=\gamma_{a}\circ J_{a}$ we denote the 1-form defined by $Z\longmapsto \gamma_{a}(J_{a}Z)$. 
Moreover, the part of $T^{H, \om}(X,Y)$ in the torsion component $\mc{X}_6$ is given by
\begin{gather*}
\frac16\sum_{cycl}(\psi_a(X)J_bY+\psi_a(J_aX)J_cY-\psi_a(Y)J_bX-\psi_a(J_aY)J_cX)\,,
\end{gather*}
where $\psi_a$ is the local 1-form defined by $\psi_a:=\gamma_c\circ J_a-\gamma_b$.  \et
 \pr
 Since the difference of $\nabla^{H, \om}_X$ and $\nabla^{Q, \om}_X$ is 1-form valued in $\fr{sp}(1)\oplus \fr{so}^*(2n)$, $\nabla^{H, \om}$ preserves $Q$ and $\om$.  Next we check which 1-form  $-\frac12\sum^3_{a=1}\tilde\gamma_a(X)J_a$ valued in $\fr{sp}(1)$ corrects $\nabla^{Q, \om}$ to preserve the local almost complex structure $J_a\in H$. To this end, we solve the following equation: 
 \begin{align*}
0&=(\nabla^{Q, \om}J_a)Y-\frac12(\tilde\gamma_aJ_a+\tilde\gamma_bJ_b+\tilde\gamma_cJ_c)J_aY+\frac12J_a(\tilde\gamma_aJ_a+\tilde\gamma_bJ_b+\tilde\gamma_cJ_c)Y\\
&=\gamma_{c}J_{b}Y-\gamma_{b}J_{c}Y-\frac12(-\tilde\gamma_a-\tilde\gamma_bJ_c+\tilde\gamma_cJ_b)Y+\frac12(-\tilde\gamma_a+\tilde\gamma_bJ_c-\tilde\gamma_cJ_b)Y\\
&=(\gamma_{c}-\tilde\gamma_c)J_{b}Y-(\gamma_{b}-\tilde\gamma_b )J_{c}Y\,,
\end{align*}
for any $X, Y\in\Gamma(TM)$. Hence we should have $\gamma_{c}=\tilde\gamma_c$, $\gamma_{b}=\tilde\gamma_b$ and this explains the term $-\frac12\sum^3_{a=1}\gamma_a(X)J_a$ in (\ref{1.7}).  \\
\noindent  Consider now some non-trivial  $\zeta\in [\E\Hh]^*$.  According to \cite[Theorem 4.12]{CGWPartI} a  component of the form $\sum^3_{a=1}\zeta\circ J_a\otimes J_a\subset [\E\Hh]^*\otimes \sp(1)$   needs to be corrected by $-4\pi_S(\omega\otimes \zeta^T)\subset [\E\Hh]^*\otimes \fr{so}^*(2n)$ in order to satisfy the condition of \cite[Theorem 5.3]{CGWPartI} for the minimal adapted connections $\nabla^{H,\omega}$ and $\nabla^{Q, \om}$.  
Thus, in the expression $-4\pi_S(\omega\otimes \zeta^T)$ we have 
\[
\zeta^{T}=\frac{1}{6}(\sum^3_{a=1}\gamma_a(J_a))^T
\]
and  the expression in (\ref{1.7}) follows by Theorems 4.12 and 5.3 in \cite{CGWPartI}. 
 The relation between $T^{H, \om}$ and $T^{Q, \om}$ is   then a direct consequence of (\ref{1.7}).

\noindent Let us now determine the $\mc{X}_6$-part of the torsion $T^{H, \om}$. 
 The torsion $T^{Q, \om}$  of $\nabla^{Q, \om}$ and the term 
\[
-\frac23(\pi_S(\omega\otimes (\sum^3_{a=1}\gamma_a(J_a))^T)(X,Y)-\pi_S(\omega\otimes (\sum^3_{a=1}\gamma_a(J_a))^T)(Y,X))
\]
cannot contribute to this  torsion type (see      \cite[Section 5]{CGWPartI}). 
Therefore, it suffices to compute the $\mc{X}_6$-component  only of the expression  $\frac12\sum^3_{a=1}(\gamma_a(Y)J_aX-\gamma_a(X)J_aY)$ appearing in $T^{H, \om}$.  To determine this precisely, we first compute
\begin{align*}
\pi_{J_a}(\gamma_a(.)J_a)(X,Y)&=\frac14(\gamma_a(X)J_aY+J_a(\gamma_a(J_aX)J_aY+\gamma_a(X)J_aJ_aY)-\gamma_a(J_aX)J_aJ_aY)\\
&=\frac14(\gamma_a(X)J_aY-\gamma_a(J_aX)Y-\gamma_a(X)J_aY+\gamma_a(J_aX)Y)=0\,,\\
\pi_{J_c}(\gamma_a(.)J_a)(X,Y)&=\frac14(\gamma_a(X)J_aY+J_c(\gamma_a(J_cX)J_aY+\gamma_a(X)J_aJ_cY)-\gamma_a(J_cX)J_aJ_cY)\\
&=\frac14(\gamma_a(X)J_aY+\gamma_b(J_cX)J_bY+\gamma_a(X)J_aY+\gamma_a(J_cX)J_bY)\\
&=\frac12(\gamma_a(X)J_aY+\gamma_a(J_cX)J_bY)\,,\\
\pi_{J_b}(\gamma_a(.)J_a)(X,Y)&=\frac14(\gamma_a(X)J_aY+J_b(\gamma_a(J_bX)J_aY+\gamma_a(X)J_aJ_bY)-\gamma_a(J_bX)J_aJ_bY)\\
&=\frac14(\gamma_a(X)J_aY-\gamma_a(J_bX)J_cY+\gamma_a(X)J_aY-\gamma_a(J_bX)J_cY)\\
&=\frac12(\gamma_a(X)J_aY-\gamma_a(J_bX)J_cY)\,.
\end{align*}
Therefore,
\begin{align*}
\pi_{J_a}(\sum^3_{a=1}\gamma_a(.)J_a)(X,Y)&=\frac12(\gamma_b(X)J_bY+\gamma_b(J_aX)J_cY+\gamma_c(X)J_cY-\gamma_c(J_aX)J_bY)\\
&=-\frac12(\psi_a(X)J_bY+\psi_a(J_aX)J_cY)\,,
\end{align*}
and the given expression of the $\mc{X}_{6}$-component of $\nabla^{H, \om}$ follows.  
 \pro
The following is   an immediate  consequence of   Lemma \ref{integrable_J_a} and Theorem \ref{them_1.7}.

\bc\label{corol_Ja_integrable} 
Let $(M, Q, \om)$ be an almost qs-H manifold and let $H=\{J_a : a=1, 2, 3\}$ be an admissible basis of $Q$.
Then, the local almost complex structure $J_a$ is integrable if and only if the following two conditions hold:
\begin{gather}
\pi_{J_{a}}(T^{Q, \om})(X, Y)=0\,,\label{1.8}\\
\psi_a(X)J_bY+\psi_a(J_aX)J_cY-\psi_a(Y)J_bX-\psi_a(J_aY)J_cX=0\,,\label{1.9}
\end{gather}
for all $X, Y\in\Gamma(TM)$, where $\psi_a$ is the local 1-form defined by $\psi_a:=\gamma_c\circ J_a-\gamma_b$. 
\ec


\section{Submanifolds of  almost quaternionic skew-Hermitian manifolds}\label{section2}
Given an almost  quaternionic skew-Hermitian manifold $(M^{4n}, Q, \om)$, one can consider several types of submanifolds. 
In this section we will undertake a systematic study of  \textsf{(almost) complex  submanifolds}, \textsf{(almost) symplectic submanifolds}, several types of \textsf{(almost) pseudo-Hermitian submanifolds} and finally    \textsf{(almost) quaternionic submanifolds} of  $(M^{4n}, Q, \om)$.

We  initiate  our study  with almost symplectic submanifolds of $(M^{4n}, Q, \om)$. 
To facilitate this, first we describe   results concerning the extrinsic geometry of almost symplectic submanifolds of \textsf{almost Fedosov manifolds}.

\subsection{Almost symplectic submanifolds}
\subsubsection{Extrinsic geometry of almost symplectic submanifolds}\label{Extrinsic_geom1}
\noindent Let us fix an almost Fedosov manifold $(M, \om, \nabla^\om)$. For our goals it is sufficient to assume that $\dim M>4$.
 Let $N\subset M$ be a submanifold  of $M$ and let us denote by  $\iota : N\to M$ the    immersion and by $\hat\om:=\iota^*\om$ the pull-back 2-form on $N$. In the sequel we will adopt the   following definition:
\bd
(1) A submanifold $N\subset M$ of an almost Fedosov manifold $(M, \om, \nabla^\om)$ is called \textsf{almost symplectic} if the \textsf{induced 2-form} $\hat\om\in\Om^2(N)$ is non-degenerate everywhere on $N$. Such a submanifold is denoted by $(N, \hat\om)$. \\
(2)  An almost symplectic submanifold  $(N, \hat\om)$ of  $(M, \om, \nabla^\om)$ whose induced 2-form  $\hat\om$ is closed is called a  \textsf{symplectic submanifold} of $M$. 
\ed

\br
 Observe that the dimension of an almost symplectic submanifold $(N, \hat\om)$ of $(M, \om, \nabla^\om)$ is necessarily even. Moreover, the induced 2-form $\hat\om$ establishes an isomorphism $TN\cong T^*N$.
\er

  Let us fix an almost symplectic submanifold $(N, \hat\om)$  of an almost Fedosov manifold $(M, \om, \nabla^\om)$.  To simplify our description, we  may regard $\iota$ as an inclusion and identify $x=\iota(x)$ and $T_{x}N=\dd\iota_{x}(T_{x}N)=\iota_{*}(T_{x}N)$. 
 Then, at any $x\in N$ we can consider the decomposition   (see also \cite{CGRS09})
  \[
   T_{x}M=T_{x}N\oplus (T_{x}N)^{\perp_{\om}}=T_{x}N\oplus \nu_{x}(N)\,,
  \]
where $\nu_{x}(N):= (T_{x}N)^{\perp_{\om}}$ denotes 
the  \textsf{$\om$-complement} of $T_{x}N$ in $ T_{x}M$,
  \[
  \nu_{x}(N):= (T_{x}N)^{\perp_{\om}}=\{u\in T_{x}M : \om(u, v)=0\, \ \text{for all} \  v\in T_{x} N\}\,.
  \]
Hence  we have a canonical choice of normal subspaces $N\ni x\longmapsto \nu_{x}(N)\subset T_{x}M$. 
We will refer to  $\nu_{x}(N)$ as the  \textsf{normal space} of $T_{x}N$ in $T_{x}M$, and  consider the \textsf{normal bundle} $\nu(N)$ of $N$, defined by
\[
\nu(N):=\bigsqcup_{x\in N}(T_{x}N)^{\perp_{\om}},\]
such that $TM|_{N} = TN\oplus \nu(N)$. 

\br
Saying that $(N, \hat\om)$  is an almost symplectic submanifold  of $(M, \om, \nabla^\om)$ is equivalent requiring $T_{x}N$ to be a symplectic subspace of the symplectic vector space $(T_{x}M, \om_x)$, at any $x\in N$. If this is the case, the $\om$-orthogonal decomposition $T_{x}M=T_{x}N\oplus\nu_{x}(N)$ implies that $\nu_{x}(N)=(T_{x}N)^{\perp_{\om}}$ is also a symplectic subspace of $(T_{x}M, \om_x)$, with linear symplectic form given by the restriction of $\om_x$. We will denote this restriction by $\tilde\om_x$, so that   $\om_x=\hat\om_x+\tilde\om_x$ at any $x\in N$.
\er

Let us now consider the natural projections  $\top : TM\to TN$ and $\perp : TM\to \nu(N)$ to the tanget bundle and normal bundle, respectively, induced by the splitting $TM|_{N}=TN\oplus\nu(N)$.   For  any $X, Y\in\Gamma(TN)$  we set $\tilde\nabla_{X}Y:=(\nabla^{\om}_{X}Y)^{\top}$ and  define the \textsf{second fundamental form} $\al$ of the almost symplectic submanifold  $(N, \hat\om)$ by the rule $\al(X, Y):=(\nabla^{\om}_{X}Y)^{\perp}$.  Note that  $\al(X, Y)\in\Gamma(\nu(N))$ for any $X, Y\in\Gamma(TN)$ and thus we have the splitting 
\begin{equation}\label{topperp}
\nabla^{\om}_{X}Y=(\nabla^{\om}_{X}Y)^{\top}+(\nabla^{\om}_{X}Y)^{\perp}=\tilde\nabla_{X}Y+\al(X, Y)\,.
\end{equation}
We will refer to this formula as the \textsf{symplectic Gauss formula}. 
The rule $\tilde\nabla_{X}Y:=(\nabla^{\om}_{X}Y)^{\top}$ defines a connection on $N$ which  preserves the almost symplectic form $\hat\om=\iota^*\om$. In particular, 
\bp\label{induced_connection} \textnormal{(\cite{CGRS09})}
 Let $(N, \hat\om)$ be an almost symplectic submanifold of an almost Fedosov manifold $(M, \om,  \nabla^\om)$.
 Then the connection $\tilde\nabla$ is an almost symplectic connection on $(N, \hat\om)$. When $\nabla^{\om}$ is torsion-free, then $\tilde\nabla$ is a symplectic connection 
 on the symplectic submanifold $(N, \hat\om)$. 
\ep
 In \cite{CGRS09} this result  is  mentioned only for the torsion-free (i.e.,  symplectic) case, without a proof. Of course,  an approach follows the  Riemannian case of  isometric immersions, by replacing  the Levi-Civita connection  by a symplectic connection. For the almost symplectic case the approach is similar and  it suffices to use an almost symplectic connection. 
 Although this is a standard fact we give the proof for self-completeness. 

  \pr
  Since $\nabla^{\om}$ is an almost symplectic connection, we have $\nabla^{\om}\om=0$ which means that
 \begin{equation}\label{om_par}
0=(\nabla^{\om}_{X}\om)(Y, Z)=X(\om(Y, Z))-\om(\nabla^{\om}_{X}Y, Z)-\om(Y, \nabla^{\om}_{X}Z)
 \end{equation}
for any $X, Y, Z\in TM$.

\noindent Suppose now that $X, Y, Z$ are tangent vector fields on $N$, i.e., $X, Y, Z\in TN$. 
Then, by (\ref{topperp}) we have $\nabla^{\om}_{X}Y=\underbrace{\tilde\nabla_{X}Y}_{\in TN}+\underbrace{\al(X, Y)}_{\in\nu(N)}$ 
and by the definition of $\nu(N)$ we see that 
\[
\om(\nabla^{\om}_{X}Y, Z)=\om(\tilde\nabla_{X}Y+\al(X, Y), Z)=\om(\tilde\nabla_{X}Y, Z)
\]
and similarly $\om(Y, \nabla^{\om}_{X}Z)=\om(Y, \tilde\nabla_{X}Z)$.  Since $\hat\om=\om|_{TN}$,  by (\ref{om_par}) we thus obtain
\[
0=X(\om(Y, Z))-\om(\tilde\nabla_{X}Y, Z)-\om(Y, \tilde\nabla_{X}Z)=X(\hat\om(Y, Z))-\hat\om(\tilde\nabla_{X}Y, Z)-\hat\om(Y, \tilde\nabla_{X}Z)=(\tilde\nabla_{X}\hat\om)(Y, Z)\,,
\]
for  any $X, Y, Z\in TN$. This proves that $\tilde\nabla$ is an almost symplectic connection on $(N^{2k}, \hat\om)$. 
  The claim for the torsion-free case follows easily; the distribution $TN\subset TM$ is involutive and hence integrable, and in particular  we have $[X, Y]^{\top}=[X, Y]$ and $[X, Y]^{\perp}=0$, for any $X, Y\in TN$.
   \pro
 
\bd
The connection $\tilde\nabla$ is called  the \textsf{induced  connection} on $(N, \hat\om)$ by the almost symplectic connection $\nabla^\om$.   
\ed

\bp\label{second_fund_sym}
 The second fundamental form $\al$ of an almost symplectic submanifold  $(N, \hat\om)$ of an almost Fedosov manifold $(M, \om, \nabla^\om)$ satisfies
\begin{equation}\label{fund_forma}
\al(X, Y)-\al(Y, X)=(T^{\om}(X, Y))^{\perp}\,,\quad X, Y\in\Gamma(TN)\,,
\end{equation}
where $T^{\om}$ is the vector-valued torsion 2-form of the almost symplectic connection $\nabla^{\om}$. Moreover,  the torsion $\tilde T$  of the induced connection $\tilde\nabla$ satisfies the relation 
\[
\tilde T(X, Y)=(T^{\om}(X, Y))^{\top}\,,\quad X, Y\in\Gamma(TN)\,.
\]
 \ep
\pr
For any $X, Y\in\Gamma(TN)$ we have $[X, Y]=[X, Y]^{\top}$. Moreover, the  decomposition $TM|_{N}=TN\oplus\nu(N)$ implies  that  
\[
 T^{\om}(X, Y)=(T^{\om}(X, Y))^{\top}+(T^{\om}(X, Y))^{\perp}
 \]
   for any $X, Y \in\Gamma(TN)$.
Thus the  identity in (\ref{fund_forma})    is consequence of the symplectic Gauss formula,  
\begin{eqnarray*}
\al(X, Y)-\al(Y, X)&=&\nabla^{\om}_{X}Y-\nabla^{\om}_{Y}X-(\nabla^{\om}_{X}Y-\nabla^{\om}_{Y}X)^{\top}\\
&=&T^{\om}(X, Y)+[X, Y]-(T^{\om}(X, Y)+[X, Y])^{\top}\\
&=&T^{\om}(X, Y)-(T^{\om}(X, Y))^{\top}=(T^{\om}(X, Y))^{\perp}\,.
\end{eqnarray*}
Then, for the torsion $\tilde T$ of $\tilde\nabla$ we obtain
\begin{eqnarray*}
\tilde T(X, Y)&=&T^{\om}(X, Y)-\big(\al(X, Y)-\al(Y, X)\big)=T^{\om}(X, Y)-(T^{\om}(X, Y))^{\perp}=(T^{\om}(X, Y))^{\top}
\end{eqnarray*}
for all $X, Y\in\Gamma(TN)$. Alternatively,  
\[
\tilde{T}(X, Y)=\tilde\nabla_{X}Y-\tilde\nabla_{Y}X-[X, Y]=(\nabla^{\om}_{X}Y)^{\top}-(\nabla^{\om}_{Y}X)^{\top}-[X, Y]^{\top}=(T^{\om}(X, Y))^{\top}
\]
 for any $X, Y\in\Gamma(TN)$. 
\pro

\bc\label{corol_funda_al}
  Let $(N, \hat\om)$ be an almost symplectic submanifold of an almost Fedosov manifold $(M, \om,  \nabla^\om)$. Then the following hold:
\begin{itemize}
\item[(1)] If  $\al$ is symmetric, then the torsion $T^{\om}$ satisfies $T^{\om}(X, Y)=\tilde{T}(X, Y)$ for any $X, Y\in\Gam(TN)$. 
\item[(2)]   If $(M, \om, \nabla^\om)$ is a Fedosov manifold, then  the second fundamental form $\al$ is a symmetric bilinear form, that is, $\al\in\Gam(S^{2}T^*N\otimes\nu(N))$. 
\end{itemize}
\ec

  Let $(N, \hat\om)$ be an almost symplectic submanifold of an almost Fedosov manifold $(M, \om,  \nabla^\om)$. 
   If the second fundamental form vanishes, i.e., $\al=0$ identically, then $(N, \hat\om)$ is called  an \textsf{almost symplectic totally geodesic submanifold} of   $ (M, \om, \nabla^{\om})$.
   In this case  we have  $T^{\om}(X, Y)=\tilde{T}(X, Y)$ for any $X, Y\in\Gam(TN)$.
 
For an almost symplectic submanifold   $(N, \hat\om)$  of an almost Fedosov manifold $(M, \om,  \nabla^\om)$ one can also introduce a connection $\nabla^{\perp}$ on the normal bundle $\nu(N)$. In particular, for any section $\xi\in\Gamma(\nu(N))$ and any vector field $X\in\Gamma(TN)$ we may define
\[
\nabla^{\perp}_{X}\xi:=(\nabla^{\om}_{X}\xi)^{\perp}\,,\quad A_{\xi}(X):=(\nabla^{\om}_{X}\xi)^{\top}\,.
\]
  Note that   $\nabla^{\perp}$ is only a connection on a vector bundle, referred to as the   \textsf{normal connection}.  However, in an analogous way one can show that
\bp The connection $\nabla^{\perp}$  preserves the almost symplectic structure $\tilde\om_x$ at the fibers $\nu_{x}(N)$ at any $x\in N$.

\ep
The endomorphism $A_{\xi}$ is called the \textsf{shape operator} of the almost symplectic submanifold $(N, \hat\om)$.  Note that the second fundamental  form and the shape operator are  related by the identity
\[
\hat\om_{x}(A_{\xi}X, Y)=\tilde\om_{x}(\al_{x}(X, Y), \xi)\,,\quad \xi\in (T_{x}N)^{\perp_{\om}}, \ X, Y\in T_{x}N\,,
\]
at any $x\in N$.

  \subsubsection{Almost symplectic submanifolds of almost qs-H  manifolds}\label{almost_symplectic_subQ1}
 \noindent  Let us now fix an almost  quaternionic skew-Hermitian manifold $(M^{4n}, Q,  \om)$. As we mentioned in the introduction, next we may assume--without loss of generality--that  $n>1$ (to keep our text simple, we will avoid repeating this condition). 
 The remainder of this subsection is devoted to almost symplectic submanifolds  $(N^{2k}, \hat\om)$ of $(M^{4n}, Q,  \om)$,  for some $k<2n$, where as usual $\hat\om$ is the induced 2-form on $N$.

 Let us first pose a result regarding such submanifolds,  based on the  classification of intrinsic torsion components of $\SO^*(2n)\Sp(1)$-structures, as described in the introduction (see Table \ref{Table1}).
 \bp
Let  $(N^{2k}, \hat{\om})$   be an almost symplectic submanifold of an almost quaternionic skew-Hermitian manifold $(M^{4n}, Q, \om)$. If $(M^{4n}, Q, \om)$ is torsion-free or has intrinsic torsion type $\mc{X}_{15}$, then 
$(N^{2k}, \hat{\om})$  is a symplectic submanifold of $M$.
\ep

\pr
When $\nabla^{Q, \om}$ is torsion-free, the scalar 2-form $\om$ is a symplectic 2-form (see Theorem \ref{connections}). Thus, for an almost symplectic  submanifold $(N^{2k}, \hat{\om})$  of a quaternionic skew-Hermitian manifold $(M^{4n}, Q, \om)$ the 2-form $\hat{\om}$ is automatically symplectic.  More in general, we know by  \cite[Theorem 1.5]{CGWPartII} that an almost quaternionic skew-Hermitian manifold
  $(M^{4n}, Q, \om)$ has symplectic scalar 2-form $\om$ if and only if has   intrinsic torsion type $\mc{X}_{15}$, see also Table \ref{Table1}.  This implies our statement. 
 \pro
  Observe that since we are dealing with genuine immersions of lower dimension, i.e., $k<2n$, the pullback of $\om$ does not detect all components of the differential $\dd\om$. Hence, it is possible to have $\dd\hat{\om}=0$ without $\om$ being closed on $M^{4n}$, and  in general a converse statement  fails.
  
Let $(N^{2k}, \hat{\om})$  be an almost symplectic submanifold of an almost qs-H manifold $(M^{4n}, Q, \om)$ and let $\nabla^{Q,  \om}$ be the almost qs-H connection on $(M^{4n}, Q, \om)$ described in Theorem \ref{connections}. 
Consider the corresponding $\om$-orthogonal decomposition with respect to the scalar 2-form $\om\in\Om^{2}(M)$, i.e.,
 \[
  T_{x}M= T_{x}N\oplus (T_{x}N)^{\perp_{\om}}=T_{x}N\oplus \nu_{x}(N)\,.
  \]
    For any  $X, Y\in\Gamma(TN)$ we can thus set $\tilde\nabla_{X}Y:=(\nabla^{Q, \om}_{X}Y)^{\top}$ and  $\al(X, Y):=(\nabla^{Q, \om}_{X}Y)^{\perp}$, 
such that 
\[
\nabla^{Q, \om}_{X}Y=\tilde\nabla_{X}Y+\al(X, Y)\,.
\]
The following result is a consequence of Proposition \ref{induced_connection} and Proposition \ref{second_fund_sym}.
 
 \bc \label{almost_sym_sub_cor}
 Let $(N^{2k}, \hat\om)$ be an almost symplectic submanifold of an almost quaternionic skew-Hermitian manifold $(M^{4n}, Q, \om)$ with $k\leq 2n$.
 Then the induced connection $\tilde\nabla$ by $\nabla^{Q, \om}$ is an almost symplectic connection on $(N, \hat\om)$. Its torsion satisfies the identity 
 \begin{equation}\label{tilde_torsion1}
 \tilde T(X, Y)=(T^{Q, \om}(X,Y))^{\top}\,,\quad X, Y\in\Gamma(TN)\,,
 \end{equation}
 where $T^{Q, \om}$ is the torsion of $\nabla^{Q, \om}$. 
 When $\nabla^{Q, \om}$ is torsion-free, i.e., $(M^{4n}, Q, \om)$ is a qs-H manifold,  then $\tilde\nabla$ is a symplectic connection 
 on the symplectic submanifold $(N, \hat\om)$.
 \ec
For the torsion-free case, we can construct a variety of (homogeneous) explicit examples of symplectic submanifolds  of qs-H manifolds,
in terms of  qs-H symmetric spaces.  We  postpone the presentation of such examples to   Section \ref{section4}, see  Proposition \ref{ExThm1} and Example \ref{Ex_SOstar}.

\subsection{Almost complex submanifolds} 

\subsubsection{Compatible almost complex submanifolds}\label{cacsub}
\noindent Let $(M^{4n}, Q,  \om)$ be an almost  quaternionic skew-Hermitian manifold.  
In this subsection we will study almost complex submanifolds of $(M^{4n}, Q, \om)$ 
which are compatible with $Q$, 
in the following sense:

\bd
 A $2k$-dimensional submanifold    $N^{2k}$  of an almost  quaternionic skew-Hermitian manifold   $(M^{4n}, Q,  \om)$ is said to be
 an  {\sf almost complex submanifold} of $M^{4n}$, when $N$ admits an almost complex structure 
 $\hat{J} : TN\to TN$ which satisfies the following property:
 \begin{itemize}
 \item[(C1)] For any $x\in N$ there exists $J\in\Gamma(Q_{x})$ such that $J(T_{x}N)=T_{x}N$ and $J|_{T_{x}N}=\hat{J}_{x}$ at any $x\in N$. 
 \end{itemize}  
Such a submanifold will be denoted by $(N^{2k}, \hat{J})$. When  $\hat{J}$ is a complex structure, then $(N, \hat{J})$ is said to be a  {\sf complex submanifold} of $(M^{4n}, Q, \om)$.
\ed
Therefore,  for an almost complex submanifold $(N^{2k}, \hat{J})$ of   $(M^{4n}, Q,  \om)$  the tangent bundle $TN$ of $N$ should be $J$-invariant and the almost complex structure   $\hat{J}$ should be the restriction on $TN$ of a section $J\in\Gamma(Q)$.
  Similarly with the quaternionic K\"ahler case (see \cite{AM00a, AM002}), the section $J\in\Gamma(Q|_{N})$ is uniquely determined by $\hat{J}$.
 In fact, for any $x\in N$,  one can  choose an admissible basis $H=\{J_a : a=1, 2, 3\}$ of $Q$
 defined in an open neighbourhood $U\subset M$ of $x$ in $M$ such that 
 \[
 J_{1}|_{T(N\cap U)}=\hat{J}\,.
 \]
 Such an admissible basis of $Q$ will be referred to as a (local) \textsf{adapted basis} of $(N^{2n}, \hat{J})$. Note that  for simplicity one may further assume that   $U\supset N$.  
 Since the open set $U$ will not play a very important role below, we will   just write  $J_{1}|_{TN}=\hat{J}$.

Let    $(N^{2k}, \hat{J})$  be an almost complex  submanifold of an almost quaternionic skew-Hermitian manifold $(M^{4n}, Q, \om)$. 
We may characterize the integrability of $\hat{J}$ in terms of the almost qs-H connection $\nabla^{Q, \om}$ of $(M^{4n}, Q, \om)$. To describe this approach, let us consider    an admissible basis  $\{J_a : a=1, 2, 3\}$ of $Q$.  Since $\nabla^{Q, \om}$ preserves $Q$, 
we may use  (\ref{quat_conn}), i.e., the identity $\nabla^{Q, \om}_{X}J_{a}=\gamma_{c}(X)J_{b}-\gamma_{b}(X)J_{c}$
  for any $X\in\Gamma(TM)$ and any cyclic permutation $(a, b, c)$ of $(1, 2, 3)$. 
  We first establish the following useful lemma:
  \bl\label{scomplex}
Let  $(N^{2k}, \hat{J})$ be an almost complex  submanifold of an almost quaternionic skew-Hermitian manifold $(M^{4n}, Q, \om)$. Let $\{J_a : a=1, 2, 3\}$ be an admissible basis of $Q$ such that  $\hat{J}=J_1|_{TN}$, i.e.,  $\{J_a : a=1, 2, 3\}$ is an adapted basis for $(N^{2k}, \hat{J})$. Then  the following are equivalent:  
\begin{itemize}
\item[(1)] $\nabla^{Q, \om}_{J_1X}J_1-J_{1}\nabla^{Q, \om}_{X}J_1=0$  for all $X\in\Gamma(TN)$.
\item[(2)]  The locally defined 1-form $\psi:=\gamma_{3}\circ J_1-\gamma_{2}$ satisfies $\psi(X)=0$ for all $X\in\Gamma(TN)$.  
\end{itemize}
\el
This correspondence  relies on the fact that $\nabla^{Q, \om}$ is an almost qs-H connection, and hence an almost quaternionic connection.  A  similar identity holds  when $(N^{2k}, \hat{J})$ is an almost complex submanifold of a quaternionic K\"ahler manifold $(M^{4n}, Q, g)$, where our connection $\nabla^{Q, \om}$ is replaced by the Levi-Civita connection $\nabla^g$,  see  \cite[Theorem 4.1]{AM00a} (or \cite[Theorem 1.1]{AM002}).  For completeness, we present the proof.
\pr
By (\ref{quat_conn}) and the relations $J_{1}J_2=J_3$, $J_1J_3=-J_2$, we see that
\begin{eqnarray*}
\nabla^{Q, \om}_{J_1X}J_1-J_{1}\nabla^{Q, \om}_{X}J_1&=&\big\{\gamma_{3}(J_1X)J_2-\gamma_{2}(J_1X)J_3\big\}-J_1\big\{\gamma_{3}(X)J_2-\gamma_{2}(X)J_3\big\}\\
&=&\gamma_{3}(J_1X)J_2-\gamma_{2}(J_1X)J_3-\gamma_{3}(X)J_3-\gamma_{2}(X)J_2\\
&=&\big\{\gamma_{3}(J_1X)-\gamma_{2}(X)\big\}J_2-\big\{\gamma_{2}(J_1X)+\gamma_{3}(X)\big\}J_3\,.
\end{eqnarray*}
However, $\gamma_{3}(J_1X)-\gamma_{2}(X)=(\gamma_{3}\circ J_1-\gamma_{2})(X)$ for all $X\in\Gamma(TN)$  and since $\hat{J}(X)=J_1(X)$, we also obtain that
\begin{eqnarray*}
\gamma_{2}(J_1X)+\gamma_{3}(X)&=&\gamma_{2}(\hat{J}X)+\gamma_{3}(X)=\gamma_{2}(\hat{J}X)-\gamma_{3}(\hat{J}^{2}X)=-(\gamma_{3}\circ \hat{J}-\gamma_{2})(\hat{J}X)\\
&=&-(\gamma_{3}\circ J_1-\gamma_{2})(\hat{J}X)\,.
\end{eqnarray*}
Therefore,  a combination of these relations gives
\begin{eqnarray}\label{super_c}
\nabla^{Q, \om}_{J_1X}J_1-J_{1}\nabla^{Q, \om}_{X}J_1&=&\big\{\gamma_{3}\circ J_1-\gamma_{2}\}(X)J_2+\big\{\gamma_{3}\circ J_1-\gamma_{2}\}(\hat{J}X)J_3\nonumber\\
&=&\psi(X)J_2+\psi(\hat{J}X)J_3\,.
\end{eqnarray}
The result now easily follows.
\pro
 
For simplicity, next we will write  $\uppsi:=\psi|_{TN}=(\gamma_{3}\circ J_1-\gamma_{2})|_{TN}$ for the locally defined 1-form associated with the adapted basis, restricted to $TN$.

\bt\label{Thm_1form}
Let  $(N^{2k}, \hat{J})$ be an almost complex  submanifold of an almost quaternionic skew-Hermitian manifold $(M^{4n}, Q, \om)$,
with  an adapted basis $\{J_a : a=1, 2, 3\}$.  
Then, for $k=1$ the almost complex structure $\hat{J}$ is automatically integrable. For $k>1$ 
 the following hold: 
\begin{itemize}
\item[(1)]  The  almost complex structure  $\hat{J}$ is integrable 
if and only if 
\[
\uppsi=0\,,\quad\text{and}\quad \pi_{J_{1}}(T^{Q, \om})(X, Y)=0\,,\quad\text{for all}\quad X, Y\in\Gamma(TN)\,. 
\]
\item[(2)]  Suppose that   $(M^{4n}, Q, \om)$ is torsion-free, i.e., a qs-H manifold. Then $\hat{J}$ is integrable 
if and only if the 1-form $\psi$ vanishes on $TN$, i.e.,  $\uppsi=0$.
\item[(3)] More in general, suppose that   $(M^{4n}, Q, \om)$ is of type $\mc{X}_{345}$, i.e.,  $Q$ is quaternionic.  Then $\hat{J}$ is integrable 
if and only if the 1-form $\psi$ vanishes on $TN$, i.e.,  $\uppsi=0$.
\end{itemize}
\et

\pr
It is well-known that  any 
almost complex structure in dimension 2 is integrable. Let us now assume that $k>1$.   By the theory of almost complex submanifolds of almost complex manifolds we know that one can identify the  Nijenhuis tensor $N_{\hat{J}}$ of   $\hat{J}=J_{1}|_{TN}$  on $N$
 with the restriction of  the Nijenhuis tensor $N_{J_{1}}$ to $TN$, that is, 
\[
N_{\hat{J}}(X, Y)=N_{J_{1}}(X, Y)\,,\quad X, Y\in\Gamma(TN)\,.
\]
Thus the vanishing of $N_{\hat{J}}$ is equivalent to the vanishing of $N_{J_{1}}|_{TN}$.
For the integrability of the local almost complex structure $J_1$ we can thus use  Corollary \ref{corol_Ja_integrable},
which gives that $J_1$ is integrable if and only if
the conditions  (\ref{1.8}) and (\ref{1.9}) hold, which now  both need to be restricted   to $TN$.
In particular, by Lemma \ref{scomplex} we can write 
\begin{eqnarray}
N_{\hat{J}}(X, Y)&=&N_{J_{1}}(X, Y)\nonumber\\
&=&\big\{(\nabla^{Q, \om}_{J_1X}J_1)Y -J_1(\nabla^{Q, \om}_{X}J_1)Y\big\}-\big\{(\nabla^{Q, \om}_{J_1Y}J_1)X-J_{1}(\nabla^{Q, \om}_{Y}J_1)X\big\}\nonumber\\
&& - T^{Q, \om}(J_1X, J_1Y) + J_1T^{Q, \om}(J_1X, Y) + J_{1}T^{Q, \om}(X, J_1Y) + T^{Q, \om}(X, Y)\nonumber\\
&\overset{(\ref{super_c})}{=}&\psi(X)J_2(Y)+\psi(\hat{J}X)J_3(Y)-\psi(Y)J_2(X)-\psi(\hat{J}Y)J_3(X)
 +4\pi_{J_{1}}(T^{Q, \om})(X, Y)\,,\nonumber\\ \label{The_N_J}
\end{eqnarray}
for all $X, Y\in\Gamma(TN)$. Hence   by Corollary \ref{corol_Ja_integrable} it follows  that the condition $N_{\hat{J}}(X, Y)=0$ for all $X, Y\in\Gamma(TN)$
is equivalent to the conditions
\[
\pi_{J_{1}}(T^{Q, \om})(X, Y)=0\,,\quad \psi_1(X)J_2Y+\psi_1(J_1X)J_3Y-\psi_1(Y)J_2X-\psi_1(J_1Y)J_3X=0
\]
for all $X, Y\in\Gamma(TN)$.\\
\noindent  We will now show that  the vanishing of the second quantity is equivalent with the condition $\psi|_{TN}=0$.
Obviously, the condition $\psi|_{TN}=0$  implies the vanishing of the latter expression. 
Conversely, fix some arbitrary $x\in N$. Since   $k>1$, we have $\dim T_{x}N\geq 4$ and   we may choose a 4-dimensional subspace ${\sf W}:={\rm span}\{u, J_1u,v, J_1v\}\subset T_{x}N$,  for some tangent vectors $u, v\in T_{x}N$. 
For $X=u$ and $Y=v$ we get
\[
\psi_1(u)J_2v+\psi_1(J_1u)J_3v-\psi_1(v)J_2u-\psi_1(J_1v)J_3u=0
\]
  and since the vectors $J_2u,J_3u,J_2v,J_3v $ are  linearly independent (as they are the $J_2$-images of the  linearly independent vectors $u, J_1u, v, J_1v$), it follows that 
  \[
  \psi_1(u)=0\,,\quad \psi_1(J_1u)=0\,,\quad \psi_1(v)=0\,,\quad \psi_1(J_1v)=0\,,
  \]
  that is, $\psi|_{{\sf W}}=0$.  Since $x\in N$  and the vectors $u,v$ were arbitrary (in particular, ${\sf W}\subset T_{x}N$ is arbitrary), every tangent vector in $T_{x}N$
can be included in such a subspace ${\sf W}$. We thus conclude that $\psi|_{TN}=0$. \\ 
\noindent The   claim  in part (2) for the case where $\nabla^{Q, \om}$ is torsionless is now direct.  Finally, for part (3) let us  recall that (see \cite{G97, CGWPartI})
\[
{\sf Ker}(\pi_H)=\cap_{a=1, 2, 3}{\sf Ker}(\pi_{J_{a}})\,,
\]
where $H=\{J_{a} : a=1, 2, 3\}$ is the local almost hypercomplex structure. Since $\mc{X}_{3457}\subset{\sf Ker}(\pi_{H})$ it follows by the previous relation  that $\pi_{J_1}(T^{Q, \om})=0$. This   completes the proof. 
\pro

\br
We remark that, in the torsion-free case, an analogous result hold 
 for almost complex submanifolds $(N^{2k}, \hat{J})$ of a quaternionic K\"ahler manifold $(M^{4n}, Q, g)$,  see for example part (1) in  \cite[Theorem 4.1]{AM00a},
 or part (1) in \cite[Theorem 1.1]{AM002}. 
Later on, this result will  be reformulated in terms of   almost pseudo-Hermitian submanifolds of almost qs-H manifolds.
\er

\subsubsection{Remarks regarding the torsion-free case and the $\mc{X}_{345}$-class}\label{t_FREE1}
 \noindent Let  $(M^{4n}, Q, \om)$ be an  almost qs-H manifold and  let $(N^{2k}, \hat{J})$ be an  almost complex submanifold with an adapted basis $\{J_a : a=1, 2, 3\}$.  In this subsection we assume that $k>1$. 
For any point $x\in N$, we denote by $T^{Q}_{x}N$   the \textsf{maximal $Q_{x}$-invariant subspace} of $T_{x}N$. Hence, $T^{Q}_{x}N\subset T_{x}N$ is maximal between all linear subspaces $\Wg\subset T_{x}N$ of $T_{x}N$ with respect to the property  $J_{x}\Wg=\Wg$ for any $J_{x}\in\Gamma(Q_x)$ and $x\in N$. Obviously, $T_{x}^{Q}N$ is $\hat{J}_{x}$-invariant, and under our assumptions we have 
   \[
 T^{Q}_{x}N=J_{2}T_{x}N\cap T_{x}N\,.
   \]
   We  mention  that  the dimension of these subspaces may not be constant.

\bp\label{lemma_TQN}
Let  $(N^{2k}, \hat{J})$ $(k>1)$  be an almost complex  submanifold of an almost quaternionic skew-Hermitian manifold $(M^{4n}, Q, \om)$. 
 Set  
 \begin{equation}\label{Vxy}
(V_{X, Y}^{\psi})_{x}:=\psi_{x}(X)Y-\psi_{x}(\hat{J}X)\hat{J}Y-\psi_{x}(Y)X+\psi_{x}(\hat{J}Y)\hat{J}X\,,
\end{equation}
for all  $X, Y\in T_{x}N$. If  $(M, Q, \om)$ is torsion-free or of type $\mc{X}_{345}$, then  both the tangent vectors $(N_{\hat{J}}(X, Y))_{x}$ and $(V_{X, Y}^{\psi})_{x}$ belong to the $Q$-invariant directions of $T_{x}N$, that is, $(N_{\hat{J}}(X, Y))_{x}\in T^{Q}_{x}N$ and $(V_{X, Y}^{\psi})_{x}\in T^{Q}_{x}N$ for all $X, Y\in T_{x}N$. 
 \ep
 \pr
Let $\{J_{a} : a=1, 2, 3\}$  be an   adapted basis for   $(N^{2k}, \hat{J})$,  such that $\hat{J}=J_1|_{TN}$. 
Then, the relation $J_{3}=-J_{2}J_1$ gives $J_{3}(X)=-J_{2}(J_{1}X)=-J_{2}(\hat{J}(X))$, for all $X\in\Gamma(TN)$. Thus,  by the proof in Theorem \ref{Thm_1form}    we obtain   
\begin{eqnarray*}
(N_{\hat{J}}(X, Y))_{x}&=&\psi_{x}(X)J_2(Y)+\psi_{x}(\hat{J}X)J_3(Y)-\psi_{x}(Y)J_2(X)-\psi_{x}(\hat{J}Y)J_3(X)+4(\pi_{J_{1}}(T^{Q, \om})(X, Y))_{x}\\
&=&J_{2}\big\{\psi_x(X)Y-\psi_x(\hat{J}X)\hat{J}Y-\psi_x(Y)X+\psi_x(\hat{J}Y)\hat{J}X\big\}+4(\pi_{J_{1}}(T^{Q, \om})(X, Y))_{x}
\end{eqnarray*}
for any $X, Y\in T_{x}N$. 
This shows that
\begin{equation}\label{our_NJ}
(N_{\hat{J}}(X, Y))_{x}=J_{2}\big((V_{X, Y}^{\psi})_{x}\big)+4(\pi_{J_{1}}(T^{Q, \om})(X, Y))_{x}\in T_{x}N\,,
\end{equation}
where $(V_{X, Y}^{\psi})_{x}$ is defined in (\ref{Vxy}).
When $M$ is torsion-free or  of type $\mc{X}_{345}$, then we have 
 \[
 \pi_{J_{1}}(T^{Q, \om})(X, Y)_{x}=0
 \]
  for all $X, Y\in T_{x}N$. 
Hence  in this case  (\ref{our_NJ})  gives $(N_{\hat{J}}(X, Y))_{x}=J_{2}\big((V_{X, Y}^{\psi})_{x}\big) \in J_{2}(T_{x}N)\cap T_{x}N$. Moreover, this  shows that $(V_{X, Y}^{\psi})_{x}\in T^{Q}_{x}N$ for any  $X, Y\in T_{x}N$ (see also the proof of \cite[Lemma 1.2] {AM002} for a similar situation).
 \pro
 
 \bc
 Let $(M^{4n}, Q, \om)$ be  an almost qs-H manifold which is either torsion-free or of type $\mc{X}_{345}$. If $(N^{2k}, \hat{J})$ $(k>1)$ is an almost complex submanifold of $M$ which is  totally-real with respect to $J_2$, then $\hat{J}$ is integrable.
 \ec
 \pr
 If $N$ is a totally-real submanifold of $M$ with respect to $J_2$, then we have $J_{2}T_{x}N\cap T_{x}N=\{0\}$, hence $T^{Q}N=\{0\}$. 
 Thus the assertion follows by Proposition \ref{lemma_TQN}.
\pro
 
\br  For the torsion-free case Proposition \ref{lemma_TQN}  has   an analogue  for almost complex submanifolds of quaternionic K\"ahler manifolds, see for example the proof of  \cite[Lemma 1.2]{AM002}. It  shows that for both these two torsion-free cases the obstruction to the integrability of the almost complex structure $\hat{J}$  lies in the $Q$-invariant  part  $T_x^{Q}N$ of the tangent space $T_xN$ for all $x\in N$.
\er
Let us now set 
 \[
  T^{\uppsi}_{x}N:=\Ker(\uppsi_x)\cap\Ker(\uppsi_x\circ \hat{J}_x)\subset T_{x}N
 \]
 for the intersection of the kernels of the local 1-forms $\uppsi_x$ and $\uppsi_x\circ \hat{J}_x$, for any $x\in N$. 
 \bl\label{3.20} 
 At any point $x\in N$ the subspace $T^{\uppsi}_{x}N$ is the maximal $\hat{J}$-invariant
 subspace of $T_{x}N$ where $\uppsi$ vanishes.  In particular, if $\uppsi_{x}\neq 0$ then $T^{\uppsi}_{x}N$ has codimension 2 in $T_{x}N$.
 \el
 \pr
 Obviously, a tangent vector  $u\in T^{\uppsi}_{x}N$ satisfies  $\uppsi(u)=0$ and $\uppsi(\hat{J}_{x}u)=0$. Applying $\hat{J}_{x}$ to these relations we get $\uppsi(\hat{J}_{x}u)=0$ and $\uppsi(\hat{J}_{x}^{2}u)=-\uppsi(u)=0$. Hence $T^{\uppsi}_{x}N$ is a $\hat{J}$-invariant subspace of $T_{x}N$  and  in particular  the largest $\hat{J}$-invariant subspace of $T_{x}N$ on which $\uppsi$ vanishes. The maximality follows by the definition of $T^{\uppsi}_{x}N$.   For the second claim notice that when the corresponding 1-form $\uppsi_x$ and $\uppsi_x\circ \hat{J}_x$ is nonzero at $x$, then each kernel in the definition of $T^{\uppsi}_{x}N$ is a hyperplane of $T_{x}N$. This gives  that  ${\rm codim}T^{\uppsi}_{x}N:=\dim T_{x}N-\dim T_{x}^{\uppsi}N\leq 2$ at any $x\in N$. We  obtain the equality  case since  $\uppsi_x$ and $\uppsi_x\circ \hat{J}_x$ are linearly independent, at any $x\in N$.  \pro
 
 \bc\label{3.21}
 Let  $(N^{2k}, \hat{J})$ $(k>1)$ be an almost complex  submanifold of an almost quaternionic skew-Hermitian manifold $(M^{4n}, Q, \om)$ which is torsion-free or more in general of type $\mc{X}_{345}$. Then the following hold:
 \begin{itemize}
 \item[(1)]  
At any point $x\in N$ the  following  conditions are equivalent:
  \[
  (N_{\hat{J}})_{x}\neq 0 \ \Longleftrightarrow \ \uppsi_{x}\neq 0   
  \ \Longleftrightarrow \ {\rm codim}\ T_{x}^{\uppsi}N=2\,.
  \]
 \item[(2)] At any point $x\in N$ the  following  conditions are equivalent:
  \[
  (N_{\hat{J}})_{x}=0  \ \Longleftrightarrow \ \uppsi_{x}= 0   
  \ \Longleftrightarrow \ {\rm codim}\ T_{x}^{\uppsi}N=0\,.
  \]
\item[(3)] Let $U\subset N$ be an open neighbourhood of $x\in N$ where   $\uppsi_x\neq 0$.  
 Then $T_{y}^{\uppsi}N\subset T^{Q}_{y}N$ for any $y\in U$. Moreover, if $\dim N=4k$ then $T^{Q}_{y}N=T_{y}N$ for any $y\in U$, and if $\dim N=4k+2$, then $T^{\uppsi}_{y}N=T^{Q}_{y}N$ for any   $y\in U$. 
 \end{itemize}
 \ec
 \pr
We obtain parts (1) and (2) by combining  Theorem \ref{Thm_1form} with  Lemma \ref{3.20}. 
 For part (3), we can choose a  vector $Z\in T_{x}N$ such that $\uppsi(Z)=1$ and $\uppsi(\hat{J}Z)=0$ (otherwise we could have $\Ker(\uppsi_x)=\Ker(\uppsi_x\circ\hat{J}_x)$). Then, for $Y\in T_{x}^{\uppsi}N$ and $X=Z$   the vector $(V^{\psi}_{X, Y})_{x}$ introduced in (\ref{Vxy}) takes the form
 \[
 (V^{\psi}_{Z, Y})_{x}=\psi_{x}(Z)Y-\psi_{x}(\hat{J}Z)\hat{J}Y-\psi_{x}(Y)Z+\psi_{x}(\hat{J}Y)\hat{J}Z=Y.
 \] 
Under our assumptions,   we have $(V^{\psi}_{Z, Y})_{x}\in T^{Q}_{x}N$ (see Proposition \ref{lemma_TQN}). Therefore, we deduce that $Y\in T^{Q}_{x}N$, that is, $T_{x}^{\uppsi}N\subset T_{x}^{Q}N$.
 The remaining  two claims follow easily by dimensional arguments. 
\pro

We will explore the 1-form $\uppsi$ further  during the next subsection, where we will assume  in addition  that  on the almost complex submanifold $(N, \hat{J})\subset (M, Q, \om)$ the 2-form $\hat\om:=\iota^*\om$ is an almost symplectic form. 
As it will become clear shortly, in this case $\hat\om$ is necessarily $\hat J$-invariant.

\subsection{Almost pseudo-Hermitian submanifolds}
\subsubsection{Compatible almost pseudo-Hermitian submanifolds}
 \noindent  Let us  consider an almost complex submanifold $(N^{2k}, \hat{J})$  of  an almost  qs-H manifold $(M^{4n}, Q, \om)$   which,  in addition, is  an almost  symplectic  submanifold of  $M^{4n}$. Hence, next we will additionally assume that  the induced 2-form
 $\hat{\om}\in\Om^2(N)$ is non-degenerate everywhere on $N$.

   \bp\label{almost_complex_sub}
   Let $(M^{4n}, Q, \om)$  be an almost qs-H manifold and let  $(N^{2k}, \hat{J})$ be an almost complex  submanifold of $M$,  which is additionally almost  symplectic.
   Then the pair $(\hat{g}, \hat{J})$  is an almost pseudo-Hermitian structure on $N^{2k}$, where $\hat{g}$ is defined by $\hat{g}(X, Y):=\hat{\om}(X, \hat{J}Y)$ for
any $X, Y\in\Gamma(TN)$.  
   \ep
   \pr
Since all our considerations are local,
 we may assume that $\hat{J}=J_1|_{TN}$ on some open neighbourhood $U\supset N$, where $J_1$
 in an element in an admissible basis $H=\{J_1, J_2, J_1J_2\}$ of $Q$, defined on  $U$. Hence $\{J_a\}$ is an adapted basis of $(N, \hat{J})$ and we have
  \begin{equation}\label{eqj1}
 \iota_{*}(\hat{J}(X))=J_{1}(\iota_{*}X)\,,\quad \forall \ X\in\Gamma(TN)\,,
\end{equation}
where $\iota : N\to M$ is the corresponding immersion. 
 Then, for any $X, Y\in\Gamma(TN)$ we see that
\begin{eqnarray*}
\hat{\om}(\hat{J}X, \hat{J}Y)&=&(\iota^*\om)(\hat{J}X, \hat{J}Y)=\om\big(\iota_{*}(\hat{J}X), \iota_{*}(\hat{J}Y)\big)\\
&\overset{(\ref{eqj1})}{=}&\om\big(J_{1}(\iota_{*}X), J_{1}(\iota_*Y)\big)\overset{(\ast)}{=}\om(\iota_*X, \iota_*Y)=\hat{\om}(X, Y)\,,
\end{eqnarray*}
where $(\ast)$ relies on the fact that  $\om$ is $Q$-Hermitian, hence $\{J_a\}$-invariant (see \cite[Proposition 2.10]{CGWPartI}).
This prove that  $\hat\om$ is $\hat{J}$-invariant. 
Now,  the tensor $\hat{g}$ is non-degenerate and symmetric. Since $\hat{\om}(X, Y)=\hat\om(\hat{J}X, \hat{J}Y)$ we see that
  \[
  \hat{g}(\hat{J}X, \hat{J}Y)=\hat{\om}(\hat{J}X, \hat{J}^2Y)=\hat{\om}(X, \hat{J}Y)=\hat{g}(X, Y)\,, \quad X, Y\in\Gamma(TN)\,.
  \]
  This shows that $(\hat{J}, \hat{g})$ is an almost pseudo-Hermitian structure on $N^{2k}$. 
     \pro

     \bc\label{clas_corol_1}
 Let $(M^{4n}, Q, \om)$  be an  almost qs-H manifold and let  $(N^{2k}, \hat{J}, \hat{\om},\hat{g})$  be an almost pseudo-Hermitian submanifold   of $M$ with $k>1$.  Then the following hold:
\begin{itemize}
\item[(1)]  If $(M^{4n}, Q, \om)$ is of type $\fr{X}_{15}$, i.e., the scalar 2-form $\om$ is symplectic, then $(N^{2k}, \hat{J}, \hat{\om},\hat{g})$ is   almost pseudo-K\"ahler.  
\item[(2)]  If  $(M^{4n}, Q, \om)$ is either torsion-free or of type $\mc{X}_{5}$, then   $(N^{2k}, \hat{J}, \hat{\om},\hat{g})$    is  pseudo-K\"ahler  if and only if the local 1-form $\psi$ vanishes on $TN$. 
 \end{itemize}
\ec
\pr
   (1)   If $\om$ is symplectic then  $\hat\om$ is also  symplectic and by   Proposition \ref{almost_complex_sub} it follows that    $(N^{2k}, \hat{J}, \hat{\om},\hat{g})$ is an  almost pseudo-K\"ahler manifold.  \\
    (2) If  $(M^{4n}, Q, \om)$ is torsion-free then $(N^{2k}, \hat{J}, \hat{\om},\hat{g})$ is immediately almost K\"ahler and the statement follows by part (2) in Theorem \ref{Thm_1form}. Similarly, since $\mc{X}_{345}\cap\mc{X}_{15}=\mc{X}_{5}$, when 
$(M^{4n}, Q, \om)$ is of   type $\mc{X}_{5}$  the claim is a consequence of part (3) in Theorem \ref{Thm_1form}. 
\pro

According to  \cite{GrHerv},  almost pseudo-K\"ahler manifolds form one of the four  pure intrinsic torsion classes of almost pseudo-Hermitian structures.  These pure classes are   labeled by $\mc{W}_1, \ldots, \mc{W}_4$; in total there are  $2^4=16$ classes of almost pseudo-Hermitian structures,  usually referred to as the \textsf{Gray-Hervella classes}.
We can generalize the previous corollary  using the pure Gray-Hervella classes. Let us recall that  the intrinsic torsion of an almost pseudo-Hermitian structure $(g, J, \om)$  consists of the Nijenhuis tensor $N_{J}$ and the 3-form $\dd\om$. The four pure classes are   characterized as follows:
 \begin{itemize}
       \item  $\mc{W}_1$ corresponds to  nearly pseudo-K\"ahler manifolds.  Here  $N_{J}$ is totally skew-symmetric and  $\pi_{J}(T^{\om})=T^{\om}$, where   $\om(T^{\om}(X, Y), Z)=\dd\om(X, Y, Z)$.
       \item $\mc{W}_2$ corresponds to almost pseudo-K\"ahler manifolds, i.e., $\dd\om=0$. Here the torsion is still defined by $N_{J}$ but its skew-symmetric part vanishes.
       \item $\mc{W}_3$ corresponds to  balanced almost pseudo-Hermitian structures. Here $N_{J}$ vanishes and $\dd\om$ is trace-free.
       \item $\mc{W}_4$ corresponds to locally conformal pseudo-K\"ahler manifolds. Here $N_{J}$ vanishes and $\dd\om$ is of pure-trace, that is, $\dd\om=\om\wedge\theta$, where $\theta$ is the Lee form.
       \end{itemize}
     Combining this intrinsic torsion characterization with the intrinsic torsion characterization of $\SO^*(2n)\Sp(1)$-structures presented in \cite{CGWPartI} (see also Table \ref{Table1}), we can prove  the following
\bt\label{Gr_H}
  Let $(M^{4n}, Q, \om)$  be an  almost qs-H manifold and let  $(N^{2k}, \hat{J}, \hat{g})$  be an almost pseudo-Hermitian submanifold   of $M$. 
  Let us denote by $T^{Q, \om}_{ab}$ the torsion component of $T^{Q, \om}$ lying in the intrinsic  torsion  part $\mc{X}_{ab}$ of  $(M^{4n}, Q, \om)$, and let   $V^{\psi}$ be the bilinear form associated to the vector defined in {\rm(\ref{Vxy})}.
 Then the following hold:
\begin{itemize}
\item[(1)]  The intrinsic torsion of   $(\hat{J}, \hat{g})$ corresponding to $\mc{W}_1$ is given by the projection $\pi_{J_1}(T^{Q, \om}_{2})$, restricted to $N$.
\item[(2)] The intrinsic torsion of  $(\hat{J}, \hat{g})$  corresponding to  $\mc{W}_2$  is given by  the expression   $\pi_{J_1}(T^{Q, \om}_{1})+J_{2}(V^{\psi})$, restricted to $N$.
\item[(3)] The intrinsic torsion of  $(\hat{J}, \hat{g})$  corresponding to  $\mc{W}_3$   is given by the projection $(\id-\pi_{J_1})(T^{Q, \om}_{23})$, restricted to $N$.
\item[(4)]  The intrinsic  torsion of  $(\hat{J}, \hat{g})$  corresponding to  $\mc{W}_4$ is given by $T^{Q, \om}_{4}\in\mc{X}_4$, restricted to $N$.
       \end{itemize}
  \et   
  \pr
 Let us recall by part (2) of \cite[Corollary 5.7]{CGWPartI} that the 3-form $\dd\om$ corresponds to  torsion  lying in the component $\mc{X}_{234}$, and this is the skew-symmetric part of $T^{Q, \om}$.  Moreover, the Nijenhuis tensors $N_{\hat{J}}$ is characterized by  (\ref{our_NJ}), and only the torsion components in $\mc{X}_{12}$ contribute.\\
 \noindent (1) Obviously,  $\mc{X}_{12}\cap\mc{X}_{234}=\mc{X}_{2}$. Therefore,  the torsion of the almost Hermitian structure $(\hat{J}, \hat{g}, \hat{\om})$  in this case comes from the 
 projection  $\pi_{J_1}(T^{Q, \om}_{2})$, restricted to $N$. \\
  \noindent (2) Here the torsion is the complementary part of Nijenhuis tensor to the torsion of the part (1), i.e., to $\mc{X}_2$. Since only the torsion components in $\mc{X}_{12}$ contribute,  we thus  get  the $\mc{X}_1$-part $T^{Q, \om}_{1}$ of $T^{Q, \om}$. \\  
    \noindent (3) Since $\mc{X}_4$ encodes the vectorial part of  $T^{Q, \om}$ (see \cite[Corollary 5.7]{CGWPartI}), the component of $\dd\om$ in $\mc{X}_{23}$ complementary to the torsion component of  the part (1) gives the torsion in this case. In particular, we have $\mc{X}_{3}\subset{\sf Ker}(\pi_{J_1})$, therefore $\pi_{J_1}(T^{Q, \om}_{23})=\pi_{J_1}(T^{Q, \om}_{2})$ and the torsion in this case is given by the projection $(\id-\pi_{J_1})(T^{Q, \om}_{23})$. \\
     \noindent (4) The  $\mc{X}_{4}$-torsion component  is the vectorial part $\dd\om$. Thus  the claim   follows.
        \pro

        Using this we obtain the following 
        \bc\label{Gr_H_corol}
        Let $(M^{4n}, Q, \om)$  be an  almost qs-H manifold and let  $(N^{2k}, \hat{J}, \hat{g})$  be an almost pseudo-Hermitian submanifold   of $M$. Then, 
                \begin{itemize}
\item[(1)]   $(N^{2k}, \hat{J}, \hat{g})$ is nearly pseudo-K\"ahler  (type $\mc{W}_1$) if and only if  we have
\[
0=\pi_{J_1}(T^{Q, \om}_{1})+J_{2}(V^{\psi})=(\id-\pi_{J_1})(T^{Q, \om}_{23})= T^{Q, \om}_{4}\quad\text{on}\quad TN\,.
\]
\item[(2)] $(N^{2k}, \hat{J}, \hat{g})$ is almost pseudo-K\"ahler (type $\mc{W}_2$) if and only if   $T_{234}^{Q, \om}|_{TN}=0$.
\item[(3)]  $(N^{2k}, \hat{J}, \hat{g})$ is balanced (type $\mc{W}_3$) if and only if  we have
\[
0=\pi_{J_1}(T^{Q, \om}_{2})=\pi_{J_1}(T^{Q, \om}_{1})+J_{2}(V^{\psi})= T^{Q, \om}_{4}\quad\text{on}\quad TN\,.
\]
\item[(4)]  $(N^{2k}, \hat{J}, \hat{g})$ is  locally conformal pseudo-K\"ahler  (type $\mc{W}_4$) if and only if  we have
\[
T_{23}^{Q, \om}|_{TN}=0=\big(\pi_{J_1}(T^{Q, \om}_{1})+J_{2}(V^{\psi})\big)|_{TN}\,.
\]
\item[(5)] The pseudo-K\"ahler  case is characterized as follows:
\begin{itemize}
\item[$\bullet$] If $k=1$ then $(N^2, \hat{J}, \hat{g})$ is automatically a  K\"ahler manifold (a Riemann surface).
\item[$\bullet$]  If $k=2$ then $(N^4, \hat{J}, \hat{g})$ is a pseudo-K\"ahler manifold if and only if $\uppsi=0$, $\pi_{J_1}(T^{Q, \om}_{1})=0$ and $T^{Q, \om}_{4}=0$ on $TN$. 
\item[$\bullet$]  If $k>2$ then $(N^{2k}, \hat{J}, \hat{g})$ is pseudo-K\"ahler manifold if and only if only the component $T^{Q, \om}_{15}$ has non-trivial restriction to $TN$, 
  $\uppsi=0$ and $\pi_{J_1}(T_1^{Q, \om})=0$ on $TN$.
\end{itemize} 
 \end{itemize}
\ec
       
       \br
  (1) For $k=1$ all the modules $\mc{W}_1, \ldots, \mc{W}_4$ are trivial, while for $k=2$  the modules $\mc{W}_1$ and $\mc{W}_3$ are trivial. Thus the generic case is for $k>2$ (see \cite{GrHerv} for the Riemannian case and observe that the results in \cite{GrHerv} are applicable to the more general case of $\U(p, q)$-structures, because   the Lie groups $\U(n)$ and $\U(p, q)$ have the same complexification).  \\
 (2)  According to Theorem \ref{Gr_H} it follows that the $\mc{X}_{5}$-component of the torsion  $T^{Q, \om}$ does not contribute to the classes $\mc{W}_{1}, \ldots, \mc{W}_4$.  \er 
 
 \bex
 Let $M$ be  the linear model $([\E\Hh], Q_0, \om_0)$ and let $a+bj$ with $a, b\in\C^{n}$ be the coordinates of the corresponding skew-Hermitian bases (see \cite[Definition 2.24]{CGWPartI}). Let $H_{0}=\{\mc{J}_{a} : a=1, 2, 3\}$ be the standard admissible basis of $Q_0$. Consider the metric $g_{0}^{2}=g_{\mc{J}_{2}}$  corresponding to $\mc{J}_2$ via the standard scalar 2-form $\om_0$.    By the proof in \cite[Proposition 2.11]{CGWPartI} we know that 
 this metric can be expressed as
 \[
 g_{0}^{2}(a+bj, c+dj)={\rm Re}(a^{t}c+b^{t}d)\,.
 \]
 Consider now the subset $ N_{p, q}\subset [\E\Hh]$ of $[\E\Hh]$ defined by
\begin{align*}
 N_{p, q}:=\{a+bj\in [\E\Hh] : {\rm Re}(a_{\al})=0={\rm Re}(b_{\al}) \ (1\leq \al\leq p), \text{and} \ &\\
 \quad\quad\quad\quad\quad\quad  {\rm Im}(a_{\beta})=0={\rm Im}(b_{\beta})  \ (p+1\leq \be\leq p+q)\}\,,
\end{align*}
 where $p+q=n$.
This is a $2n$-dimensional submanifold of $[\E\Hh]$ which is  $\mc{J}_{2}$-invariant. Moreover, since $g_0^2$ is a non-degenerate Hermitian product, the 2-form $\hat\om:= \om_{0}|_{TN_{p, q}}$ is non-degenerate as well. Thus $(N_{p, q}, \hat{J}:=\mc{J}_2|_{TN_{p, q}}, \hat{g}:=g_{0}^{2}|_{TN_{p, q}})$ is a pseudo-K\"ahler submanifold of $([\E\Hh], Q_0, \om_0)$ of signature $(2q, 2p)$. 
Later,  in Section \ref{section4}, we will   construct examples in the homogeneous setting, see   Proposition \ref{ExThm1}. 
 \eex

\subsubsection{Integrability conditions in terms of the induced connection}
\noindent Let $(M^{4n}, Q, \om)$  be an  almost qs-H manifold and let  $(N^{2k}, \hat{J}, \hat{\om},\hat{g})$  be an almost pseudo-Hermitian submanifold  of $M$,  as  described  in Proposition \ref{almost_complex_sub}. Based on the results of Section \ref{almost_symplectic_subQ1},  at any point $x\in N$ one can consider     the $\om$-orthogonal decomposition   
  \[
  T_{x}M= T_{x}N\oplus (T_{x}N)^{\perp_{\om}}=T_{x}N\oplus \nu_{x}(N)\,.
  \]
  \bl
  The normal bundle $\nu(N)$ is $J_1$-invariant, that is, $J_1(\nu(N)) \subset \nu(N)$.
  \el
  \pr
  Let $u\in\nu_{x}(N)$ be some non-zero vector in $\nu_{x}(N)=(T_{x}N)^{\perp_{\om}}$. Then $\om_{x}(u, X)=0$ for all $X\in T_{x}N$. 
   We also have $\om_{x}(u, J_1X)=0$ for all $X\in T_{x}N$. However, $\om$ is $Q$-Hermitian, therefore it satisfies $\om_{x}(J_{1}U, V)+\om_{x}(U, J_1V)=0$ for all $x\in M$ and $U, V\in T_{x}M$.  For $U=u\in\nu_{x}(N)\subset T_{x}M$, $V=X\in T_{x}N\subset T_{x}M$, and $x\in N$  this gives
  $\om_{x}(J_{1}u, X)+\om_{x}(u, J_1X)=0$. Since  $\om_{x}(u, J_1X)=0$, we thus  get $\om_{x}(J_{1}u, X)=0$ for all $X\in T_{x}N$. This shows that   $J_{1}u\in\nu_{x}(N)$ for all $x\in N$, hence proves our claim. 
  \pro
Let us now characterize the integrability of $\hat{J}$ using the induced connection 
  \begin{equation}\label{induced_Q_con}
  \tilde\nabla_{X}Y:=(\nabla^{Q, \om}_{X}Y)^{\top}=\nabla^{Q, \om}_{X}Y-\al(X, Y)\,,\quad   X, Y\in\Gamma(TN )\,,
  \end{equation}
  where $\al$ is the second fundamental form. Recall that $\al(X, Y)\in\Gamma(\nu(N))$.

  \bp\label{almost_symp_conq1}
  Let $(M^{4n}, Q, \om)$  be an  almost qs-H manifold and let  $(N^{2k}, \hat{J}, \hat{\om},\hat{g})$  be an almost pseudo-Hermitian submanifold   of $M$, as  described  in Proposition \ref{almost_complex_sub}. Then  the induced connection  $\tilde\nabla$  is an almost symplectic connection,    $\tilde\nabla\tilde\om=0$, whose torsion $\tilde T$ satisfies the relation {\rm(\ref{tilde_torsion1})}. 
  If $\nabla^{Q, \om}$ is torsion-free, then $\tilde\nabla$ is also torsion-free and the second fundamental form $\al$ is symmetric.
   \ep
  \pr
  This is a consequence of Corollaries   \ref{corol_funda_al} and \ref{almost_sym_sub_cor}. 
    \pro

  \bt\label{main_application_acs}
  Let $(M^{4n}, Q, \om)$  be an  almost qs-H manifold and let $(N^{2k}, \hat{J}, \hat{\om},\hat{g})$ be  an  almost pseudo-Hermitian submanifold  of $M^{4n}$, with $k>1$.  Suppose that  $\pi_{J_{1}}(T^{Q, \om})(X, Y)=0$ for all $X, Y\in\Gamma(TN)$, and moreover that one of the following equivalent conditions holds:
 \begin{itemize}
  \item[(1)] $(\tilde\nabla_{X}\hat{J})Y=0$ for all $X, Y\in\Gamma(TN)$;
  \item[(2)] For any  $X, Y\in\Gamma(TN)$ we have 
   \begin{equation}\label{integrability_hatJ2}
   \big(\gamma_{3}(X)J_2(Y)-\gamma_{2}(X)J_{3}(Y)\big)^{\top}=0\,. 
   \end{equation}
  \item[(3)] The connection $\tilde\nabla$ is metric with respect to $\hat{g}$, i.e., $\tilde\nabla\hat g=0$.
  \end{itemize}  
  Then  the almost complex structure $\hat{J}$ is integrable.
\et
\pr
In terms of the induced connection $\tilde\nabla$  on $N$ we have the relation
\begin{equation}\label{tilde_N_J_1}
N_{\hat{J}}(X, Y)=\big\{(\tilde\nabla_{\hat{J}X}\hat{J})Y -\hat{J}(\tilde\nabla_{X}\hat{J})Y\big\}-\big\{(\tilde\nabla_{\hat{J}Y}\hat{J})X-\hat{J}(\tilde\nabla_{Y}\hat{J})X\big\}+4\pi_{\hat{J}}(\tilde T)(X, Y)
\end{equation}
for all $X, Y\in\Gamma(TN)$.
Moreover, by Corollary \ref{almost_sym_sub_cor}
 we know that the torsions $\tilde{T}$ and $T^{Q, \om}$ are related by (\ref{tilde_torsion1}), that is 
$ \tilde T(X, Y)=(T^{Q, \om}(X,Y))^{\top}$, for all $X, Y\in\Gamma(TN)$.  Hence if $\pi_{J_{1}}(T^{Q, \om})(X, Y)=0$  then we will have $\pi_{\hat{J}}(\tilde T)(X, Y)=0$ for all $X, Y\in\Gamma(TN)$. Moreover, if $(\tilde\nabla_{X}\hat{J})Y=0$ for all $X, Y\in\Gamma(TN)$, then 
by (\ref{tilde_N_J_1}) we   obviously get $N_{\hat{J}}(X, Y)=0$ for all $X, Y\in\Gamma(TN)$, that is, $\hat{J}$ is integrable.\\
\noindent It now remains to prove the equivalence of the three given conditions.
Recall that  we may assume that $\hat{J}=J_1|_{TN}$, where $\{J_a : a=1, 2, 3\}$ is an adapted basis of $(N^{2k}, \hat{J})$ (with respect to $Q$).  Moreover, we have  $\al(X, Y)\in\Gamma(\nu(N))$, for any $X, Y\in\Gamma(TN)$. Hence by applying  the Gauss formula  (\ref{induced_Q_con})   we compute  \begin{eqnarray}\label{gen_induced_form}
  (\tilde\nabla_{X}\hat{J})Y&=&\tilde\nabla_{X}\hat{J}Y-\hat{J}(\tilde\nabla_{X}Y)\nonumber\\
  &=&\nabla^{Q, \om}_{X}J_1Y-\al(X, \hat{J}Y)-J_1(\nabla^{Q, \om}_{X}Y)+J_{1}(\al(X, Y))\nonumber\\
  &=&(\nabla^{Q, \om}_{X}J_1)Y-\al(X, \hat{J}Y)+J_1\big(\al(X, Y)\big)\,,
  \end{eqnarray}
for all $X, Y\in\Gamma(TN)$. 
 Therefore, one has $\tilde\nabla\hat J=0$ if and only if 
\begin{equation}\label{integrability_hatJ}
  (\nabla^{Q, \om}_{X}J_1)Y=\al(X, \hat{J}Y)-J_1(\al(X, Y))\,,\quad X, Y\in\Gamma(TN)\,.
  \end{equation}
Recall now by (\ref{quat_conn}) that we can write
 \begin{equation}\label{quat_connJ1}
 (\nabla^{Q, \om}_{X}J_1)Y=\gamma_{3}(X)J_{2}Y-\gamma_{2}(X)J_{3}Y\,.
 \end{equation}
  By restricting both sides of  (\ref{integrability_hatJ}) to $TN$  and using (\ref{quat_connJ1}) we obtain
 \[
\Big(\gamma_{3}(X)J_2(Y)-\gamma_{2}(X)J_{3}(Y)\Big)^{\top}=\Big(\al(X, \hat{J}Y)-J_1\big(\al(X, Y)\big)\Big)^{\top}=0\,, 
\]
for any $X, Y\in\Gamma(TN)$.   
This gives the equivalence between (1) and (2), namely,
\[
 (\tilde\nabla_{X}\hat{J})Y=0 \quad \Longleftrightarrow \quad  \big(\gamma_{3}(X)J_2(Y)-\gamma_{2}(X)J_{3}(Y)\big)^{\top}=0 
\]
for all $X, Y\in\Gamma(TN)$.  Let us finally prove that  the condition in (1) is equivalent with the condition in (3) and the equivalence between the conditions in (2) and (3)  follows.
Indeed, if $\tilde\nabla\hat J=0$, then together with the relation $\tilde\nabla\hat\om=0$ described in  Proposition \ref{almost_symp_conq1}, we obtain that $\tilde\nabla\hat g=0$, hence $\tilde\nabla$ is metric with respect to $\hat{g}$. Similarly, if $\tilde\nabla\hat g=0$ then together with the relation $\tilde\nabla\hat\om=0$ one can show that $\tilde\nabla\hat J=0$. This proves our theorem.
\pro

Next we will adopt the following definition (motivated by an analogous definition given in \cite{Sekigawa} for almost complex submanifolds of  almost Hermitian manifolds).
\bd
An almost pseudo-Hermitian submanifold   $(N^{2k}, \hat{J}, \hat{\om},\hat{g})$  
of an almost qs-H manifold $(M^{4n}, Q, \om)$, is said to be an \textsf{$\al$-submanifold}
if the second fundamental form $\al$ satisfies the identity 
\[
\al(\hat{J}X, Y)=\al(X, \hat{J}Y)=J_1\al(X, Y)\,,\quad X, Y\in\Gamma(TN)\,.
\]
\ed

The following is now a corollary of Theorem \ref{main_application_acs}.
\bc Let $(M^{4n}, Q, \om)$  be an  almost qs-H manifold and let $(N^{2k}, \hat{J}, \hat{\om},\hat{g})$ be  an  almost pseudo-Hermitian submanifold  of $M^{4n}$, for some $k>1$. 
\begin{itemize}
  \item[(1)] If $N^{2k}$ is totally geodesic or  an $\al$-submanifold of $M^{4n}$ and $\pi_{J_{1}}(T^{Q, \om})(X, Y)=0$, for all $X, Y\in\Gamma(TN)$,    then $\hat{J}$ is integrable if and only if  $\nabla^{Q, \om}J_1=0$.
\item[(2)]  If  $(M^{4n}, Q, \om)$ is torsion-free 
and the condition in {\rm(\ref{integrability_hatJ2})} is satisfied, then $(N^{2k}, \hat{J}, \hat{\om},\hat{g})$ is pseudo-K\"ahler.  
\end{itemize}
\ec

  \pr
 The first assertion is a consequence of  the relation (\ref{integrability_hatJ})   in combination with part  (1) of Theorem \ref{Thm_1form}. For the second one, note that  if   $\nabla^{Q, \om}$ is torsion-free,  then  so is $\tilde\nabla$, see Corollary \ref{almost_sym_sub_cor}. Hence, when $\nabla^{Q, \om}$ is torsion-free and  the integrability condition   (\ref{integrability_hatJ}) is satisfied, then by Theorem \ref{main_application_acs}  it follows that the  connection $\tilde\nabla$ should coincide with the Levi-Civita connection $\nabla^{\hat{g}}$ and the relation $\nabla^{\hat{g}}\hat{J}=0$ implies that $\hat{J}$ is integrable. 
\pro

The conditions (1)-(3)  in Theorem \ref{main_application_acs} admit another equivalent interpretation, which is  described as follows:
\bt\label{totally_complex}
Let $(M^{4n}, Q, \om)$  be an  almost qs-H manifold and let $(N^{2k}, \hat{J}, \hat{\om},\hat{g})$ be  an  almost pseudo-Hermitian submanifold  of $M^{4n}$, for some $k>1$. Then   $\tilde\nabla\hat{J}|_{x}=0$  if and only if one of the following   conditions holds at $x\in N$:
\begin{itemize}
\item[{\rm(I1)}] $\gamma_{2}|_{T_xN}=\gamma_{3}|_{T_xN}=0$;
\item[{\rm (I2)}] $J_{2}T_{x}N\perp_{\om} T_{x}N$.
\end{itemize}
\et

\smallskip
Obviously, the condition {\rm (I2)} resembles the usual  \textsf{totally-complex}  condition in the theory of almost Hermitian submanifolds of K\"ahler manifolds, where the orthogonal complement is considered with respect to a Riemannian metric. 
 \pr
 First observe that whenever one of the conditions {\rm(I1)} or {\rm(I2)} holds then the relation in  {\rm(\ref{integrability_hatJ2})}
is satisfied.  For {\rm(I1)}  this is obvious. If the condition {\rm(I2)} holds,  then   $(J_{3}(Y))^{\top}=-(J_{2}(J_{1}Y))^{\top}=0$ and $(J_2(Y))^{\top}=0$.\\
Next  we will show that  at any point $x\in N$ where $\tilde\nabla\hat{J}|_{x}=0$, at least one of the conditions {\rm (I1)}, {\rm (I2)} should be satisfied.
First note that  based on (\ref{quat_connJ1}), the relation (\ref{gen_induced_form}) for $Y$ replaced by $\hat{J}Y$ gives 
\begin{eqnarray*}
 (\tilde\nabla_{X}\hat{J})\hat{J}Y&=&(\nabla^{Q, \om}_{X}J_1)J_1Y-\al(X, \hat{J}^2Y)+J_1^2(\al(X, Y))\\
 &=&(\gamma_{3}(X)J_{2}(\hat{J}Y)-\gamma_{2}(X)J_{3}(\hat{J}Y))^{\top}\\
 &=&-(\gamma_{3}(X)J_3Y+\gamma_{2}(X)J_2(Y))^{\top}\,,
  \end{eqnarray*}
 for all  $X, Y\in\Gamma(TN)$.
 Hence,  in combination with the result of  Theorem \ref{main_application_acs}, we have proved that
 \begin{equation}\label{numbered1}
   (\tilde\nabla_{X}\hat{J})Y=(\gamma_{3}(X)J_2(Y)-\gamma_{2}(X)J_{3}Y)^{\top}\,,
   \end{equation}
   and
   \begin{equation}\label{numbered2}
    (\tilde\nabla_{X}\hat{J})\hat{J}Y=-(\gamma_{3}(X)J_3Y+\gamma_{2}(X)J_2(Y))^{\top}\,.
 \end{equation}
 When $\tilde\nabla\hat{J}|_{x}=0$, then  in the left-hand-side of these relations we get zero. \\
 \noindent  Suppose that the condition ${\rm (I1)}$ is not satisfied. We will show that ${\rm (I2)}$ is then true. Indeed,  by (\ref{numbered1})  we have  $\gamma_{2}(X)(J_3Y)^{\top}= \gamma_{3}(X)(J_2Y)^{\top}$, for all $X, Y\in\Gamma(TN)$. By substituting this equality into the $\gamma_2(X)$-multiple of (\ref{numbered2}) we then obtain  \[
 (\gamma_{2}^{2}(X)+\gamma^{2}_{3}(X))(J_{2}Y)^{\top}=0\,,
 \]
 for all $X, Y\in\Gamma(TN)$, and  our claim follows immediately.  
Similarly, if {\rm (I2)} does not hold, then  the expressions $(J_2(Y))^{\top}$ and  $(J_3(Y))^{\top}=J_1(J_2(Y))^{\top}$ for some  $Y\in\Gamma(TN)$ are linearly independent. Therefore,  we immediately see that  {\rm (I1)}  should hold.
 \pro

\br
The results established above provide an ``almost symplectic analogue'' of those obtained in the proof of Theorem 5.2 in \cite[p.~28]{AM00a} (see also Theorem 1.8 in \cite[p.~876]{AM002}),  for almost Hermitian submanifolds of quaternionic K\"ahler manifolds.
In fact, the authors in \cite{AM00a, AM002}  show that under an assumption on the reduced scalar curvature of $(M, Q, g)$, 
the conditions  {\rm(I1)} and {\rm(I2)} are equivalent. 
However, there is no counterpart to scalar curvature in our framework.
\er

\subsubsection{On the non-integrability of the almost pseudo-Hermitian structure}\label{t_FREE2}
 \noindent  Let  us fix  an almost pseudo-Hermitian  submanifold  $(N^{2k}, \hat{J}, \hat{g}=\hat{\om}\circ\hat{J})$ $(k>1)$ of an almost qs-H manifold $(M^{4n}, Q, \om)$, as described in Proposition \ref{almost_complex_sub}  and  denote by $\{J_a : a=1, 2, 3\}$ an adapted basis of $(N, \hat{J})$.  
   Let us consider the 1-form $\uppsi=\psi|_{TN}=(\gamma_{3}\circ J_1-\gamma_{2})|_{TN}$
  and the subspace  $T^{\uppsi}_{x}N:=\Ker(\uppsi_x)\cap\Ker(\uppsi_x\circ \hat{J}_x)$, introduced in Section \ref{t_FREE1}. 
Let   $\Uppsi\in\Gamma(TN)$  the (local) dual vector field of the 1-form 
  $\uppsi=\psi|_{TN}$ with respect to the almost symplectic 2-form $\hat{\om}$, that is 
  \[
  \hat{\om}(Z, \Uppsi)=\uppsi(Z)\,,\quad\text{for all}\quad Z\in \Gamma(TN)\,.
  \]
  Obviously, we have $\uppsi_x\neq 0$ if and only if $\Uppsi_{x}\neq 0$. 
As above, let us denote by
\[
\Ug^{\perp_{\hat{\om}}}=\{u\in T_{x}N : \hat\om(u, v)=0 \ \text{for all} \ v\in\Ug\}
\]
the  $\hat{\om}$-complement  of a vector subspace $\Ug\subset T_{x}N$ in the symplectic vector space $(T_{x}N, \hat{\om})$.
Recall that we always have $\dim\Ug+\dim\Ug^{\perp_{\hat{\om}}}=\dim T_{x}N$ but, unlike orthogonal complements,  it is not necessary true that $\Ug\cap \Ug^{\perp_{\hat{\om}}}=\{0\}$. When this happens, then $\Ug$ is  a symplectic subspace of $T_{x}N$, i.e., $\hat{\om}$ restricts to a non-degenerate 2-form on $\Ug$.  In  this case we have $T_{x}N=\Ug\oplus \Ug^{\perp_{\hat{\om}}}$, thus $\Ug^{\perp_{\hat{\om}}}$ is a symplectic subspace of $T_{x}N$ as well.  
 
 \bt\label{2.36}
 Let $(M^{4n}, Q, \om)$ be an  almost qs-H manifold and let  $(N^{2k}, \hat{J}, \hat{g}=\hat{\om}\circ\hat{J})$ $(k>1)$ be an almost pseudo-Hermitian  submanifold of $M$ with $\uppsi_{x}\neq 0$ at any $x\in N$ and $\pi_{J_{1}}(T^{Q, \om}_{1})(X, Y)=0$ for all $X, Y\in\Gamma(TN)$.  Then the following hold: 
 \begin{itemize}
 \item[(1)] If $\dim N=4k$, then $T_{x}N=T^{Q}_{x}N$ at any $x\in N$  and the signature of $\hat{g}$ is $(2k, 2k)$.
 \item[(2)] If $\dim N=4k+2$, then $T^{\uppsi}_{x}N=T^{Q}_{x}N$ at any $x\in N$ and the signature of $\hat{g}$ is $(2(k+1), 2k)$. Moreover, we have an $\hat\om$-orthogonal decomposition
 \begin{equation}\label{TxpsiN_split}
 T_{x}N=T_{x}^{\uppsi}N\oplus (T_{x}^{\uppsi}N)^{\perp_{\hat{\om}}}\,,\quad  (T_{x}^{\uppsi}N)^{\perp_{\hat{\om}}}={\rm span}\{\Uppsi_x, \hat{J}_{x}\Uppsi_x\}\,,
 \end{equation}
  and both the subspaces $T_{x}^{\uppsi}N=T^{Q}_{x}N$ and  $(T_{x}^{\uppsi}N)^{\perp_{\hat{\om}}}={\rm span}\{\Uppsi_x, \hat{J}_{x}\Uppsi_x\}$ are symplectic subspaces of $(T_{x}N, \hat{\om}_x)$.   In particular, in this case the (local) vector $\Uppsi_x$ is non-isotropic at any $x\in N$.
 \end{itemize}
 \et
 \pr
 By our assumptions and part (2) in Theorem 
\ref{Gr_H}  it follows that the intrinsic torsion of the pair $(\hat{J}, \hat{g})$ corresponding to $\mc{W}_2$, is given by  the restriction of   $J_{2}(V^{\psi})$   to $N$.  The given relations $T_{x}N=T^{Q}_{x}N$ in (1) and  $T^{\uppsi}_{x}N=T^{Q}_{x}N$ in (2) follow similarly with  Corollary   \ref{3.21}.   For the case with $\dim N=4k$,
the proof regarding the signature is analogous to those presented in  \cite[Proposition 2.11]{CGWPartI}. When $\dim N=4k+2$, there are two cases that we need to consider,  depending on the  character of vector  $Z\in T_{x}N$,  used in the proof of Corollary \ref{3.21}.\\
\noindent{\bf Case 1:}  We can choose $Z$ such that  $\uppsi(Z)=1$, with $Z=c\thinspace\hat{J}\Uppsi$  for some non-zero scalar $c$. 
This gives the relation
\[
1=\uppsi(Z)=\om(c\hat{J}_x\Uppsi_x, \Uppsi_x)=-c\hat{g}_{x}(\Uppsi_{x}, \Uppsi_{x})\,,
\]
thus $\hat{g}_{x}(\Uppsi_x, \Uppsi_x)\neq 0$, that is, $|\Uppsi_x|^2_{\hat g}\neq 0$.   Therefore, in this case the vector $\Uppsi_x$ is non-isotropic. 
Set  $\Vg_{x}:=T_{x}N$, $\Ug_{x}:=T_{x}^{\uppsi}N$ and $\Wg_{x}:={\rm span}\{\Uppsi_x, \hat{J}_{x}\Uppsi_x\}$.
  By definition we have
  \begin{eqnarray*}
  \Ug_{x}=T_{x}^{\uppsi}N&=&\{u\in \Vg_x : \uppsi_x(u)=0 \ \text{and} \ \uppsi_x(\hat{J}_{x}u)=0\}\\
  &=&\{u\in \Vg_x : \hat{\om}_{x}(u, \Uppsi_x)=0 \ \text{and} \ \hat{\om}_{x}(\hat{J}_{x}u, \Uppsi_x)=0\}\\
  &=&\{u\in \Vg_x : \hat{\om}_{x}(u, \Uppsi_x)=0 \ \text{and} \ \hat{\om}_{x}(u, \hat{J}_{x}\Uppsi_x)=0\}
    =\Wg_x^{\perp_{\hat{\om}}}\,,
  \end{eqnarray*}
where it was used   that the pair $(\hat{J}, \hat{g}=\hat{\om}\circ\hat{J})$ is an almost pseudo-Hermitian structure.   Thus, $\Ug_{x}=\Wg^{\perp_{\hat{\om}}}$, which gives that $\Ug_{x}^{\perp_{\hat{\om}}}=\Wg$. The condition $\Ug_{x}\cap\Wg_{x}=\{0\}$ follows since $|\Uppsi_x|^2_{\hat g}\neq 0$.  Indeed, a general vector $v\in\Ug_x\cap\Wg_x$ is written as   $v=\al\Uppsi_x+\beta\hat J_x\Uppsi_x$, 
 for some real scalars $\al,\be$, subject to the conditions
\[
0=\uppsi_x(v)=\hat\om_x(v, \Uppsi_x)=\hat\om_x(\al\Uppsi_x+\be\hat J_x\Uppsi_x, \Uppsi_x,)
=\al\hat\om_{x}(\Uppsi_x, \Uppsi_x)+ \be\,\hat\om_x(\hat J_x\Uppsi_x, \Uppsi_x,)\,,
\]
and
\[
0=\uppsi_x(\hat J_x v)=\hat\om_x(\hat J_x v, \Uppsi_x,)
=\hat\om_x(\al\hat{J}_{x}\Uppsi_x- \be\Uppsi_x, \Uppsi_x, )= -\be\,\hat\om_x(\Uppsi_x,\Uppsi_x)+\al\,\hat\om_x(\hat J_x\Uppsi_x, \Uppsi_x,)\,.
\]
By skew-symmetry we have $\hat\om_x(\Uppsi_x,\Uppsi_x)=0$, thus these equations reduce to  
\[
-\be\,\hat\omega_x(\Uppsi_x,\hat J_x\Uppsi_x)=0\,,\quad\text{and}\quad
-\al\,\hat\omega_x(\Uppsi_x,\hat J_x\Uppsi_x)=0\,,
\]
respectively. 
Since $|\Uppsi_x|^2_{\hat g}\neq 0$ we thus get $\hat\om_x(\Uppsi_x,\hat J_x\Uppsi_x)\neq0$, which gives that 
  $\al=\be=0$, that is, $\Ug_x\cap\Wg_x=\{0\}$. Since  we also have $\dim\Vg_x=\dim\Ug_{x}+\dim\Wg_x=\dim\Ug_x+2$ we finally get   the direct sum $\Vg_x=\Ug_x\oplus\Wg_x$, at any $x\in N$.
This proves the $\hat\om$-orthogonal splitting  given in (\ref{TxpsiN_split}). The claim for the signature follows for the $4k$-dimensional part $T^{Q}_{x}N$ as in the proof of \cite[Proposition 2.11]{CGWPartI}, while for the 2-dimensional part $\Wg_{x}$ we get signature  $(2, 0)$, and this gives the final statement.\\
 \noindent{\bf Case 2:}  The second possible case is when  $\mc{Q}_x:={\rm span}\{Z, \hat{J}_xZ, \Uppsi_x, \hat{J}_x\Uppsi_x\}$ is a 4-dimensional subspace of $T_{x}N$, where  $Z\in T_{x}N$  is such that $\uppsi(Z)=1$.  We will show that this leads to a contradiction.
Since $\hat\om$ is $\hat{J}$-Hermitian we see that
 \[
 \hat\om(\hat{J}Z,  \hat{J}\Uppsi)=\hat\om(Z, \Uppsi)=\uppsi(Z)=1\,.
  \]
Thus,  we conclude that the restriction of $\hat\om$  to $\mc{Q}_{x}$ is non-degenerate for any $x\in N$. Moreover, by the definition of $T^{\uppsi}_{x}N$ and since $T_x^{\uppsi}N=T^{Q}_{x}N$,  it  follows that $\hat\om(\Uppsi_x, \hat{J}_x\Uppsi_x)=0$, at any $x\in N$.  Now,  the vectors  $\Uppsi_x, \hat{J}_x\Uppsi_x$ belong to $T^{Q}_{x}N$ and $J_{2}\Uppsi_x\in T^{Q}_{x}N$ (by the definition of $T^{Q}_{x}N$), hence we deduce that   $J_{2}\Uppsi_x, J_{3}\Uppsi_x$ also belong to $T^{Q}_{x}N$. Moreover we see that $\hat\om_{x}(J_{2}\Uppsi_x, J_{3}\Uppsi_x)=0$. Next we can find scalars $\al,\be\in\R$ such that the linear combination $\Xg:=\al\Uppsi_x+\be \hat{J}\Uppsi_x + J_{2}\Uppsi_x\in T_{x}N$ is a vector of the $\hat\om$-orthogonal complement $\mc{Q}_{x}^{\perp_{\hat{\om}}}$ of $\mc{Q}_{x}\subset T_{x}N$. Note that
  \[
  \hat{J}_{x}\Xg=\al\hat{J}\Uppsi_x-\be\Uppsi_x+J_3\Uppsi_x\,.
  \]
 Thus we compute   $\hat\om(\Xg, \hat{J}\Xg)=\hat\om(\al\Uppsi_x+\be \hat{J}\Uppsi_x + J_{2}\Uppsi_x, \al\hat{J}\Uppsi_x-\be\Uppsi_x+J_3\Uppsi_x)=0$. Using  the complement $\mc{Q}_{x}^{\perp_{\hat{\om}}}$  we can now proceed with dimensional reduction, where $\Uppsi_x$ is replaced by $\Xg$. This finally ends up on a 2-dimensional subspace where $\hat\om$ is degenerate.  A contradiction.
 \pro
 
 \br
Regarding the first part of Theorem \ref{2.36}, where $\dim N=4k$,  observe that  $\Uppsi_x$ 
can be an isotropic tangent vector with respect to $\hat{g}$.
\er

\subsection{Almost quaternionic submanifolds} \label{almostqsub}
Let us finally describe  almost quaternionic submanifolds $N\subset M$ of (torsion-free) quaternionic skew-Hermitian manifolds $(M^{4n}, Q, \om)$.   
We  are thus interested in    immersed submanifolds $N^{4k}\subset M^{4n}$ whose dimension is divisible by $4$.  
Let us recall the   definition of almost quaternionic submanifolds   (see also \cite{A68, BCU81, Tasaki86,  AM00a, Ianus}).

 \bd   Let $(M^{4n}, Q)$ be an almost quaternionic manifold. A $4k$-dimensional submanifold $N^{4k}\subset M^{4n}$ of  $(M, Q)$, where $k<n$,  is called an \textsf{almost quaternionic submanifold} 
if its  tangent spaces are  $Q$-invariant, i.e., 
\[
J_{x}T_xN=T_{x}N\,, \quad \forall x\in N, J_{x}\in Q_{x}\,.
\] 
\ed
Let us recall a classical result concerning  almost quaternionic submanifolds of a quaternionic manifold $(M^{4n}, Q)$.
\bt\label{AMprop} \textnormal{(\cite{AM93, AM00a}, see also \cite[Theorem 5.3]{Pet})} 
Let $N^{4k}$ be an almost quaternionic submanifold of a quaternionic
manifold $(M^{4n}, Q)$.  Then,   $N^{4k}$ admits a quaternionic structure $\hat{Q}\subset\Ed(TN)$ given by restriction, i.e., $\hat{Q}:=Q|_{TN}$. Moreover,  let $\nabla$ be any  quaternionic connection  on $M$.  Then the restriction $\hat{\nabla}=\nabla|_{TN}$ is a quaternionic connection on $N$   and $(N, \hat{Q})$  is totally geodesic with respect to $\nabla$.
  \et

  Based on  this theorem, one usually  refers to an almost quaternionic  submanifold of a quaternionic manifold  as a \textsf{quaternionic submanifold}, see  \cite{AM00a}.
  Note that Gray \cite{Gray}, and independently  Alekseevsky \cite{A68}, 
had previously   proved Theorem \ref{AMprop}   for the special case where $M^{4n}$ is   a quaternionic K\"ahler manifold. 
As we will see below, for the case of   almost   quaternionic submanifolds of quaternionic skew-Hermitian manifolds $(M^{4n}, Q, \om)$ there is 
a similar characterization. 


 .

\bp\label{prop1}
Let $(M^{4n}, Q, \om)$   be a quaternionic skew-Hermitian manifold and let
$N^{4k}\subset M^{4n}$ be an almost quaternionic submanifold of $(M^{4n}, Q)$ for some $k<n$.   
Then 
the induced 2-form $\hat{\om}$ is $\hat{Q}$-Hermitian, where $\hat{Q}$ is the quaternionic structure  on $N^{4k}$ defined by restriction, i.e., $\hat{Q}=Q\big|_{TN}$. Moreover, $\dd\hat\om=0$.
\ep
\pr
According to Theorem \ref{AMprop},  the almost  quaternionic structure 
$\hat{Q}$ induced on $N$ by the restriction of $Q$  on the fibers of the tangent bundle $TN$, 
   is a quaternionic structure.   
Let $\{J_a : a=1, 2, 3\}$ be an admissible basis of $Q$. Since   $\hat{Q}=Q\big|_{TN}$   
and $TN$ is $Q$-invariant, 
it follows that $\{J_{a}|_{TN} : a=1, 2, 3\}$ is an admissible basis of $\hat{Q}$.
In particular,  if $\hat{H}=\{I_a\}$ is an admissible basis of $\hat{Q}$, then we may assume that  there exists an admissible basis $H=\{J_a\}$   of $Q$ such that 
\begin{equation}\label{commJa}
\iota_{*}(I_{a}X)=J_{a}(\iota_{*}X)\,, \quad a=1, 2, 3\,,
\end{equation}
for any $X\in\Gamma(TN)$, where $\iota_{*}$  is the differential of the immersion $\iota : N\hookrightarrow M$ (see for example \cite[p.~31]{AM93}). 
 Set $I_{a}=J_{a}|_{TN}$ for all $a=1, 2, 3$. Then, for the 2-form  $\hat\om=\iota^*\om$ we see that
\begin{eqnarray*}
\hat\om(I_aX, I_aY)&=&(\iota^*\om)(I_aX, I_aY)\\
&=&\om\big(\iota_{*}(I_{a}X), \iota_{*}(I_{a}Y)\big)\\
&\overset{(\ref{commJa})}{=}&\om\big(J_{a}(\iota_{*}X), J_{a}(\iota_*Y)\big)\\
&\overset{(\ast)}{=}&\om(\iota_*X, \iota_*Y)=\hat\om(X, Y)\,,
\end{eqnarray*}
for any $X, Y\in\Gamma(TN)$ and $a=1, 2, 3$.
Note here that $(\ast)$ relies on the fact that   $\om$ is $Q$-Hermitian. This proves the first claim.
The second claim about the closedness of $\hat\om$  follows   immediately, since $(M, Q, \om)$ is torsion-free, hence $\dd\om=0$.   
\pro
Based on this proposition   we  obtain the following 
\bt\label{theorem1}
Suppose that the  2-form $\hat\om$ is non-degenerate  and hence a symplectic form on $N^{4k}$.  Then, a quaternionic submanifold $N^{4k}\subset (M^{4n}, Q, \om)$ as described in Proposition {\rm \ref{prop1}}
admits a torsion-free $\SO^*(2k)\Sp(1)$-structure defined by the pair $(\hat{Q}, \hat\om)$. In particular, $(N^{4k}, \hat Q, \hat\om)$ is a quaternionic skew-Hermitian manifold
whose unique minimal adapted connection $\nabla^{\hat Q, \hat\om}$ satisfies
  $\nabla^{\hat{Q}, \hat\om}=\nabla^{Q, \om}|_{TN}=\nabla^{Q, \vol}|_{TN}$. 
\et
\pr
When the induced 2-form $\hat\om$ is non-degenerate, the quaternionic submanifold $N^{4k}\subset (M^{4n}, Q)$  described in Proposition \ref{prop1} 
admits a $\SO^*(2k)\Sp(1)$-structure $(\hat{Q}, \hat{\om})$, with $\hat{Q}$ quaternionic and $\hat{\om}$ symplectic.  Now,  by  assumption $(M, Q, \om)$ is a qs-H manifold, that is, the almost qs-H connection $\nabla^{Q, \om}$ is torsion-free, and hence a qs-H connection.
Moreover, by Theorem \ref{AMprop} it follows that the submanifold $N$ is    totally geodesic with respect to $\nabla^{Q,  \om}$, that is (see also \cite[Prop.8.2, p.~55]{Kob2} and moreover \cite[p.~33]{AM00a} for a similar case)
\[
\nabla_{X}^{Q, \om}Y\in\Gamma(TN)\,,\quad\text{for all}\quad X, Y\in\Gamma(TN)\,.
\]
This fact implies that the restricted connection $\nabla^{Q, \om}|_{TN}$ is torsion free and  hence, by uniqueness, it coincides
with the  minimal adapted connection $\nabla^{\hat{Q}, \hat{\om}}$ associated to the pair $(\hat{Q}, \hat\om)$ by Theorem \ref{connections} (the uniqueness is explained in \cite[Cor.~8.7, p.~59]{Kob2}). 
The final equality in  
\[
 \nabla^{\hat{Q}, \hat\om}=\nabla^{Q, \om}|_{TN}=\nabla^{Q, \vol}|_{TN}
\]
 follows by Theorem \ref{connections}
and the fact that $T^{Q, \om}=0$.
 This proves our claim. 
\pro
Motivated by  this result, it is natural to refer to a  quaternionic submanifold $(N^{4k}, \hat{Q}=Q|_{TN})\subset (M^{4n}, Q, \om)$  of a qs-H manifold  $(M^{4n}, Q, \om)$, which has    non-degenerate induced 2-form $\hat{\om}\in\Om^2(N)$, as a \textsf{quaternionic skew-Hermitian submanifold} of   $(M^{4n}, Q, \om)$.
We present a  homogeneous example of such  submanifolds  in the next section, see Example \ref{Ex_SOstar}.


\section{Examples}\label{section4}
In this section we will present examples based on semisimple qs-H symmetric spaces $(M=\Gg/\Lg, Q, \om)$.
Such symmetric spaces have been classified (up to covering) in \cite[Theorem 6.2]{CGWPartI} and next we will illustrate some of the results presented in Section \ref{section2}
by constructing submanifolds of such qs-H manifolds.

\subsection{Examples of submanifolds of the symmetric space $\Sl(n+1,\mathbb{H})/(\Gl(1,\mathbb{H})\Sl(n,\mathbb{H}))$}
Let us consider the 
 qs-H symmetric space  
 \[
 M=\Gg/\Lg=\Sl(n+1,\mathbb{H})/(\Gl(1,\mathbb{H})\Sl(n,\mathbb{H}))
 \]
 or real dimension $8n$, where    $\Gl(1,\mathbb{H})\Sl(n,\mathbb{H}):=\Ss(\Gl(1,\mathbb{H})\times \Gl(n,\mathbb{H}))$. 

Let  us denote by $\fr{g}=\Lie(\Gg)=\fr{sl}(n+1, \Hn)$  and $\fr{l}=\Lie(\Lg)=\fr{sp}(1)\oplus\R\oplus\fr{sl}(n, \Hn)$ the Lie algebras of $\Gg, \Lg$, respectively,  
and   let $\fr{g}=\fr{l}\oplus\fr{m}$  be the corresponding canonical decomposition. Then we have $[\fr{l}, \fr{m}]\subset\fr{m}$ and $[\fr{m}, \fr{m}]=\fr{l}$.  As usual, we will identify $\fr{m}$ with the tangent space $T_{o}M$ of $M$ at the origin $o=e\Lg\in M$.
By  \cite[Theorem 6.2]{CGWPartI}  we know that   $M$ admits a (unique) $\Sl(n+1, \Hn)$-invariant paracomplex structure $I : TM\to TM$.  This corresponds  to an ${\rm ad}_{\fr{l}}$-invariant endomorphism $I_{o} : \fr{m}\to\fr{m}$ defined by $I_{o}:={\rm ad}_{\fr{m}}(Z_0)$, where  $Z_{0}$  generates the center $ Z(\fr{l})\cong\R$  of $\fr{l}$.  Then we can define a scalar 2-form with respect to the quaternionic structure on $M$ induced by the $\fr{sp}(1)$-part in $\fr{l}$, by the relation $B_{\fr{m}}(X, Y)=\om_{o}(X, IY)$ for any $X, Y\in\fr{m}$. Here $B_{\fr{m}}$ is the restriction of the Killing form of $\fr{g}$ to $\fr{m}$.

 We can construct different  submanifolds of $M$ by considering the tower of Lie groups
 \begin{equation}\label{inclusion_RCH}
 \Sl(n+1, \R)\subset \Sl(n+1, \C)\subset \Sl(n+1, \Hn)
 \end{equation}
induced by the embeddings $\R\hookrightarrow \C\hookrightarrow \Hn$.
Both $\Sl(n+1, \R)$ and $\Sl(n+1, \C)$ are closed subgroups of $\Gg=\Sl(n+1, \Hn)$ and hence 
each inclusion in (\ref{inclusion_RCH}) induces  a homogeneous space of the form
\[
N_{\K}:=\Sl(n+1, \K)/(\Sl(n+1, \K)\cap\Lg)=\Sl(n+1, \K)/\Lg_{\K}, \quad \Lg_{\K}:=\Sl(n+1, \K)\cap\Lg\,,
\]
with $\K=\R, \C$, respectively, and $\Lg=\Gl(1,\Hn)\Sl(n,\Hn)$ is the stabilizer of $M$.
Obviously, 
\begin{eqnarray*}
\Lg_{\K}=\Sl(n+1, \K)\cap\Lg&=&\Sl(n+1, \K)\cap\Big(\Ss(\Gl(1,\mathbb{H})\times \Gl(n,\mathbb{H}))\Big)\\
&=&\Ss(\Gl(1, \K)\times\Gl(n, \K))=\Gl(1, \K)\Sl(n, \K)\,.
\end{eqnarray*}
Thus   the isotropy group $\Lg_{\K}$ has the form $\Lg_{\K}=\Sl(n+1, \K)\cap\Lg=\Gl(1, \K)\Sl(n, \K)$ and
 for $N_{\K}$ we obtain  the homogeneous  presentation  $N_{\K}=\Sl(n+1, \K)/\Gl(1, \K)\Sl(n, \K)$.  
This  is a symmetric space   for both $\R, \C$ (see for example \cite[p.~412]{Gil}). 
The corresponding canonical decomposition  is given by $\fr{g}_{\K}=\fr{sl}(n+1, \K)=\fr{l}_{\K}\oplus\fr{m}_{\K}$, 
where
 \begin{eqnarray*}
 \fr{l}_{\K}&=&\Lie(\Lg_{\K})=\left\{\begin{pmatrix} a & 0 \\ 0 & A\end{pmatrix} : a\in\K, A\in\fr{gl}(n, \K), {\rm Re}(a+{\rm tr}(A))=0\right\}\cong \fr{gl}(1, \K)\oplus\fr{sl}(n, \K)\,,\\
 \fr{m}_{\K}&=&\left\{\begin{pmatrix} 0 & u^{t} \\ w & 0 \end{pmatrix} : u, w\in\K^n\right\}\,,
\end{eqnarray*}
with  $\fr{m}_{\K}\cong\fr{g}_{\K}\cap\fr{m}$. 
Then it is easy to see that $[\fr{l}_{\K}, \fr{m}_{\K}]\subset\fr{m}_{\K}$ and $[\fr{m}_{\K}, \fr{m}_{\K}]=\fr{l}_{\K}$. Moreover, we get $[[\fr{m}_{\K}, \fr{m}_{\K}], \fr{m}_{\K}]\subset\fr{m}_{\K}$, hence $N_{\K}$ is a totally geodesic submanifold of $M$, see \cite{Kob2}.

\bp \label{ExThm1}  
 For $\K=\R$, the manifold $N_{\R}$ is a symplectic submanifold of $M$.   For $\K=\C$, the manifold $N_{\C}$ is a pseudo-K\"ahler  submanifold of $M$.
\ep
\pr  
For the first case set $\Gg_{\R}:=\Sl(n+1, \R)$ such that 
\[
N_{\R}= \Gg_{\R}/\Lg_{\R}=\Sl(n+1, \R)/(\Gl(1, \R)\Sl(n, \R))\,,\quad\text{with}\quad \dim_{\R}N_{\R}=2n\,,
\]
where $\Lg_{\R}=\Gl(1, \R)\Sl(n, \R)$. 
Let  $\iota_{\R} : N_{\R}\to M$ be the  immersion.   We know that  the center $Z(\fr{l})$  of $\fr{l}$   is generated by $Z_{0}$ and that the reductive complement $\fr{m}_{\R}$ is $I_{o}$-invariant. Now,  the Killing form on $\fr{sl}(n+1, \R)$ is the restriction of the Killing form $B$ of $\fr{sl}(n+1, \Hn)$ and the pullback $\om_{\R}:=\iota_{\R}^*\om$ is the  restriction $\om_{o}|_{TN_{\R}}$. Since we have $B_{\fr{m}_{\R}}(X, Y)=\om_{\R}(X, I_{o}Y)$ for all $X, Y\in\fr{m}_{\R}$,  it follows that $\om_{\R}$ is 
 a (closed) non-degenerate  $\Gg_{\R}$-invariant 2-form on $N_{\R}$. Thus $(N_{\R}, \om_{\R})$  is   a homogeneous symplectic submanifold of $M$. \\
\noindent For the second case,  set $\Gg_{\C}:=\Sl(n+1, \C)$ such that 
\[
N_{\C}= \Gg_{\C}/\Lg_{\C}=\Sl(n+1, \C)/(\Gl(1, \C)\Sl(n, \C))\,,\quad\text{with}\quad \dim_{\R}N_{\C}=4n\,,
\]
where $\Lg_{\C}=\Gl(1, \C)\Sl(n, \C)$. Denote by $\iota_{\C} : N_{\C}\to M$ the corresponding immersion. 
Using the same approach as for the real case,  one can show that the isotropy action of the real part $\R_{+}\subset \C^{\times}$  induces the  $\Gg_{\C}$-invariant symplectic 2-form  $\om_{\C}:=\iota_{\C}^*\om\in\Om^{2}(N_{\C})$. Hence  the pair $(N_{\C}, \om_{\C})$ is another symplectic submanifold of $M$. However,    in this case one can say more. The   $\U(1)$-part in $\Gl(1, \C)$, generates a $\G_{\C}$-invariant complex structure $\hat{J} : TN_{\C}\to TN_{\C}$ on $N_{\C}$ and it follows that  $\om_{\C}$ is $\hat{J}$-invariant.  This is because 
$\fr{u}(1)=\Lie(\U(1))\subset\fr{sp}(1)\cong Q_{o}$, hence there is an adapted basis $\{J_{a} : a=1, 2, 3\}$ of $Q$ defined on an open neighbourhood $U\subset M$, such that  $\hat{J}=J_1|_{T(N_{\C}\cap U)}$.
Then,  by  Proposition \ref{almost_complex_sub} we can show that 
 the pair $(\hat{J}, \hat{g}=\om_{\C}\circ\hat{J})$ is a $\Gg_{\C}$-invariant pseudo-Hermitian structure on $N_{\C}$, and obviously this structure is pseudo-K\"ahler.
\pro

    \br
  We should mention that not every ``reduction'' of the qs-H symmetric space $M=\Gg/\Lg=\Sl(n+1,\mathbb{H})/(\Gl(1,\mathbb{H})\Sl(n,\mathbb{H}))$   to a real subgroup  of $\Gg$ yields a symplectic submanifold of $M$. 
  Consider for instance a  maximal compact subgroup $\Kg$ of $\Sl(n+1, \Hn)$.
  Then $\Kg$ is isomorphic to $\Sp(n+1)\subset\Sl(n+1, \Hn)$,
  and the intersection $\Sp(n+1)\cap \Lg$ 
  is isomorphic to $\Sp(1)\times\Sp(n)$. This means that one can view  the quaternionic projective  $\Hn{\sf P}^{n}=\Sp(n+1)/\Sp(1)\times\Sp(n)$ as a submanifold  of $M$. However,   $\om$ does not restrict to a ($\Sp(n+1)$-invariant)  symplectic form on $\Hn{\sf P}^{n}$, since 
 the Lie algebra $\fr{sp}(1)\oplus\fr{sp}(n)$ of the isotropy group $\Sp(1)\times\Sp(n)$ has trivial center.
 In fact,   the second Betti number of $\Hn{\sf P}^{n}$ is zero, hence   this space   admits no symplectic form.
   \er

\subsection{Examples of submanifolds of the symmetric space $\SU(2+p,q)/(\SU(2)\SU(p,q)\U(1))$}
Consider the $4(p+q)$-dimensional qs-H symmetric space 
  \[
  M=\Gg/\Lg=\SU(2+p, q)/\Ss(\U(2)\times\U(p, q))\cong \SU(2+p,q)/(\SU(2)\SU(p,q)\U(1))\,.
  \]
  Let $\fr{g}=\fr{l}\oplus\fr{m}$ be the  corresponding canonical decomposition
   and let us identify $\fr{m}$ with the tangent space $T_{o}M$ of $M$ at the identity coset $o=e\Lg\in M$.
Obviously, $\fr{l}=\fr{su}(2)\oplus\fr{su}(p, q)\oplus \R Z_{0}\cong\fr{u}(2)\oplus\fr{su}(p, q)$, where $Z_0$ is the generator of the 1-dimensional center in $\fr{l}$, and $\fr{m}\cong \C^2\otimes(\C^{p+q})^{*}$ with
  \[
  \fr{m}= \left\{\begin{pmatrix} 0 & X \\ -\bar{X}^{t} & 0\end{pmatrix} : X\in{\sf Hom}_{\C}(\C^{p+q}, \C^{2})\right\}\,.
  \]
 Note that $\fr{m}$ is an irreducible $\fr{l}$-module satisfying $[\fr{m}, \fr{m}]=\fr{l}$. 
 By  \cite[Example 3.13]{CCG} 
 we know that the isotropy action of $\U(1)\subset\Lg$ gives rise to a  $\Gg$-invariant complex structure $I : TM\to TM$ with $I\notin\Gamma(Q)$, where $Q$ is the quaternionic structure on $M$ induced by the $\SU(2)$-part in $\Lg$. This corresponds to the ${\rm ad}_{\fr{l}}$-invariant endomorphism  $I_o:={\rm ad}(Z_0)$. Moreover, $M$ admits a $\SU(2+p, q)$-invariant scalar 2-form $\om$, defined via the relation
 $B_{\fr{m}}(X, Y)=\om_{o}(X, I_{o}Y)$ for all $X, Y\in\fr{m}$, where $B_{\fr{m}}$ is the restriction of the Killing form of $\fr{g}=\fr{su}(p+2, q)$ to $\fr{m}$.

  \bex \label{Lan1}
The Lie group $\hat{\Gg}:=\SO(p+2, q)$ is a closed subgroup of $\Gg=\SU(2+p, q)$ and we 
see that the intersection $\SO(p+2, q)\cap\Lg$ is isomorphic to the group $\hat{\Lg}:=\SO(2)\times\SO(p, q)\times\Z_2$.   At the level of Lie algebras we obtain  $\hat{\fr{g}}\cap\fr{l}=\Lie(\hat{\Lg})=\hat{\fr{l}}=\fr{so}(2)\oplus\fr{so}(p, q)$, where the $\fr{so}(2)$-component lies in $\fr{su}(2)$.  This yields  the homogeneous space 
     \[
   N=\hat{\Gg}/\hat{\Lg}=\SO(p+2, q)/(\SO(2)\times\SO(p, q)\times\Z_2)\,,
   \]
which is a  $2(p+q)$-dimensional totally geodesic submanifold of $M$. 
 The isotropy action of the  $\SO(2)$-component in $\hat{\Lg}$ gives rise to a $\hat{\Gg}$-invariant complex structure $\hat{J} : TN\to TN$
  on $N$. Since $\fr{so}(2)\subset\fr{su}(2)\cong\fr{sp}(1)\cong Q_{o}$ we may assume that $\hat{J}=J_{1}|_{TN}$, where $\{J_{a} : a=1, 2, 3\}$ is an admissible basis of $Q$.  
   Let $\hat{\fr{m}}$ be the reductive complement of $\hat{\fr{l}}$ in $\hat{\fr{g}}$ with respect to the Killing form   of $\hat{\fr{g}}$.     By construction we have that $\hat{\fr{m}}=\fr{m}\cap\hat{\fr{g}}$ and   $\hat{\fr{m}}$ is an  irreducible $\hat{\fr{l}}$-module. 
  By invariance, we also get  $\hat{J}_{\hat{o}}(\hat{\fr{m}})=\hat{\fr{m}}$, where $\hat{J}_{\hat{o}} : \hat{\fr{m}}\to\hat{\fr{m}}$ is the induced ${\rm ad}_{\hat{\fr{l}}}$-invariant endomorphism on $\hat{\fr{m}}=T_{\hat{o}}N$  by $\hat{J}$, with $\hat{o}=e\hat{\Lg}\in N$.   This implies that the tangent bundle $TN$ is $\hat{J}$-invariant
  and the pair $(N, \hat{J})$ is thus a homogeneous complex submanifold of $(M, Q, \om)$. \\
\noindent Note that   $\hat{\fr{m}}$ is {\it not} preserved by the $\U(1)$-factor in $\Lg$ and one can show  that 
\begin{equation}\label{tot_real1_ex}
I_{o}(\hat{\fr{m}})\cap\hat{\fr{m}}=0\,, \quad \fr{m}=\hat{\fr{m}}\oplus I_{o}(\hat{\fr{m}})\,,\quad B_{\fr{m}}(\hat{\fr{m}}, I_{o}(\hat{\fr{m}}))=0\,.
\end{equation}
 This means that  $N^{2(p+q)}$ is a totally-real submanifold of $(M^{4(p+q)}, I)$.  Let  $\hat\om_{\hat{o}} : \hat{\fr{m}}\times\hat{\fr{m}}\to\R$  be the    ${\rm ad}_{\hat{\fr{l}}}$-invariant skew-symmetric bilinear form on $\hat{\fr{m}}=T_{\hat{o}}N$ corresponding to the 2-form $\hat\om:=\iota^*\om\in\Om^2(N)$, where  $\iota : N\to M$ is the immersion.  
 Then,  by (\ref{tot_real1_ex}) it follows that $\hat\om_{\hat{o}}=0$,
 \[
 \hat\om_{\hat{o}}(X, Y)=\om_{o}(X, Y)=-B_{\fr{m}}(X, I_{o}Y)=0\,,\quad X, Y\in\hat{\fr{m}}\,.
 \]
Moreover, we have $\dim N=\frac{1}{2}\dim M$, thus $N$ is   a (complex) Lagrangian submanifold of 
 $(M^{4(p+q)}, \om)$.
\eex

 \subsection{Examples of submanifolds of the symmetric space  $\SO^*(2n+2)/\SO^*(2n)\U(1)$}\label{Examples}  
Let us finally consider the $4n$-dimensional qs-H symmetric space
\[
M^{4n}=\Gg/\Lg=\SO^*(2(n+1))/\SO^*(2n)\U(1)\,.
\]
  Let $\fr{g}=\fr{l}\oplus\fr{m}$ be the corresponding canonical decomposition,
 with $[\fr{m}, \fr{m}]=\fr{l}$ and $B(\fr{l}, \fr{m})=0$, where $B$ is the Killing form of  the Lie algebra $\fr{g}=\fr{so}^{*}(2(n+1))$ of $\Gg=\SO^*(2(n+1))$.
 We identify $\fr{m}$ with the tangent space $T_{o}\Gg/\Lg$  of $M$ at the identity coset $o=e\Lg$,
 and the  isotropy representation $\chi : \Lg\to\Aut(\fr{m})$ of $M=\Gg/\Lg$  with the restriction of the standard representation $[\E\Hh]$ of $\SO^*(2n)\Sp(1)$  to $\SO^*(2n)\U(1)$. Note that $\fr{m}$ is irreducible as an $\fr{l}$-module.    The $\Gg$-invariant quaternionic skew-Hermitian structure $(Q, \omega)$ has been  explicitly described in  \cite[Example~3.10]{CCG} (see also \cite[Section~6]{CGWPartI}).
Let us recall   that $M$ admits an  invariant complex structure $J$ which comes from a global section $M\to \mathscr{Z}$, where $\mathscr{Z}$ is the twistor space of $M$. 
This corresponds to an $\Ad(\Lg)$-invariant endomorphism  $J_{o} : \fr{m}\to\fr{m}$, defined by  $J_{o}={\rm ad}(Z_0)$ with $Z_{0}\in\fr{u}(1)$. The $\Gg$-invariant scalar 2-form $\om\in\Om^{2}(M)$ is defined in a similar manner with the previous qs-H symmetric spaces, by combining $J_{o}$ with the restriction of the  Killing form of $\fr{g}$ to $\fr{m}$.
  
 \bex\label{Ex_SOstar}
For $ k<n$ there is a standard  block diagonal embedding of the  Lie group $\hat{\Gg}:=\SO^*(2(k+1))$ into $\Gg=\SO^*(2(n+1))$, where the inclusion $\SO^*(2(k+1))\subset\SO^*(2(n+1))$ is given by
\[
\SO^*(2(k+1))\ni A\longmapsto  \begin{pmatrix} A & 0 \\ 0 & \Id_{2(n-k)}\end{pmatrix}\in\SO^*(2(n+1))\,.
\]
Obviously,  $\hat{\Gg}$ is a closed subgroup of $\Gg$ and we see that  $\hat{\Gg}\cap \Lg=\SO^*(2k)\U(1)$.
It turns out that for $k<n$ 
the homogeneous space
\[
N^{4k}=\hat{\Gg}/\hat{\Gg}\cap \Lg=\SO^*(2k+2)/\SO^*(2k)\U(1)
\]
is a totally geodesic  submanifold of $M^{4n}$.    Note that the intersection $\fr{m}\cap\fr{so}^*(2(k+1))$ is $Q$-invariant, hence we obtain a $\hat{\Gg}$-invariant quaternionic structure  $\hat{Q}$ on $N^{4k}$ by restriction, i.e.,   $\hat{Q}=Q|_{TN}$.
Moreover,  the induced  2-form $\hat\om\in\Om^{2}(N)$  is a $\hat{\Gg}$-invariant symplectic 2-form, which is $\hat{Q}$-Hermitian, i.e., a $\hat{\Gg}$-invariant scalar 2-form. 
To prove its non-degeneracy, observe that $\omega$ 
is block-diagonal with non-degenerate blocks, so its restriction to $\fr{m}\cap\fr{so}^*(2(k+1))$   is again non-degenerate. Thus,  for all $k<n$ the homogeneous space $N^{4k}$ is a quaternionic skew-Hermitian submanifold of $M^{4n}$. 
\eex

 \end{document}